\documentclass[10pt,a4paper]{article}
\usepackage{amsmath,amssymb,amsthm,graphics,graphicx,epic,eepic,multicol,ascmac}

\usepackage[all]{xy}

\numberwithin{equation}{section}
\theoremstyle{theorem}
\newtheorem{thm}{Theorem}[section]
\newtheorem{prop}[thm]{Proposition}
\newtheorem{lem}[thm]{Lemma}
\newtheorem{rem}[thm]{Remark}

\newtheorem{ex}[thm]{Example}

\theoremstyle{definition}
\newtheorem{defn}[thm]{Definition}

\def\Lm{\Lambda}
\def\nd{\noindent}
\def\ovl{\overline}

\begin{document}

\title{An algorithm for Berenstein-Kazhdan decoration functions and trails for minuscule representations}

\author{Yuki Kanakubo\thanks{Faculty of Pure and Applied Sciences, University of Tsukuba,
1-1-1 Tennodai, Tsukuba, Ibaraki 305-8577,
Japan: {y-kanakubo@math.tsukuba.ac.jp}.}, Gleb Koshevoy\thanks{Institute of Information Transmission Problems Russian Academy of Sciences, National Research University Higher School of Economics, Russian Federation :
{koshevoy@cemi.rssi.ru}.} and Toshiki Nakashima\thanks{Division of Mathematics, 
Sophia University, Kioicho 7-1, Chiyoda-ku, Tokyo 102-8554,
Japan: {toshiki@sophia.ac.jp}.} } 
\date{}
\maketitle
\begin{abstract}
For a simply connected connected simple algebraic group $G$,
a cell $B_{w_0}^-=B^-\cap U\ovl{w_0}U$ is a geometric crystal with a positive structure
$\theta_{\textbf{i}}^-:(\mathbb{C}^{\times})^{l(w_0)}\rightarrow B_{w_0}^-$.
Applying the tropicalization functor to a rational function
$\Phi^h_{BK}=\sum_{i\in I}\Delta_{w_0\Lambda_i,s_i\Lambda_i}$ called the half decoration
on $B_{w_0}^-$, one can realize the crystal $B(\infty)$ in $\mathbb{Z}^{l(w_0)}$.
By computing $\Phi^h_{BK}$, we get an explicit form of $B(\infty)$ in $\mathbb{Z}^{l(w_0)}$.
In this paper,
we give an algorithm to compute 
$\Delta_{w_0\Lambda_i,s_i\Lambda_i}\circ \theta_{\textbf{i}}^-$ explicitly for $i\in I$ such that
$V(\Lambda_i)$ is a minuscule representation of $\mathfrak{g}={\rm Lie}(G)$. In particular, the algorithm works
for all $i\in I$ if $\mathfrak{g}$ is of type ${\rm A}_n$. 
The algorithm computes a directed graph $DG$, called a {\it decoration graph}, whose vertices are labelled by all monomials in 
$\Delta_{w_0\Lambda_i,s_i\Lambda_i}\circ \theta_{\textbf{i}}^-(t_1,\cdots,t_{l(w_0)})$.
The decoration graph has some properties similar to crystal graphs of minuscule representations.
We also verify that the algorithm works in
some other cases, for example, 
the case $\mathfrak{g}$ is of type ${\rm G}_2$
though
$V(\Lambda_i)$ is non-minuscule.

\end{abstract}

\vspace{-10pt}

\section{Introduction}

The notion of `geometric crystals' was introduced in \cite{BK0}
as a geometric analog of Kashiwara's crystals, which are
defined as irreducible algebraic varieties $X$ equipped with certain
$\mathbb{C}^{\times}$-actions and rational functions, which
correspond to
Kashiwara operators and $\varepsilon$-functions, weight functions.
If there is a birational map $\theta:T'\rightarrow X$ called a positive structure
with an algebraic torus $T'$ then one can obtain a crystal via
a tropicalization functor ${\rm Trop}$ (see subsections \ref{GC-trop},\ \ref{Trop-geom}).
By this functor, $T'$ corresponds to ${\rm Trop}(T')=X_*(T')={\rm Hom}(\mathbb{C}^{\times},T')$.

We can define geometric crystal structures on varieties related to a reductive group $G$.
In \cite{BK0, BK}, it is shown
$B_{w_0}^-=B^-\cap U\ovl{w_0}U$ and $T\cdot B_{w_0}^-$ have geometric crystal structures,
where $B$, $B^-$ are opposite Borel subgroups, $U\subset B$ is a unipotent radical, $T=B\cap B^-$ is a maximal torus
and $\ovl{w_0}\in {\rm Norm}_G(T)$ is a representative of the longest element $w_0$ in Weyl group $W={\rm Norm}_G(T)/T$. 
Defining a positive structure $\theta_{\textbf{i}}: T\times (\mathbb{C}^{\times})^{l(w_0)}\rightarrow T\cdot B_{w_0}^-$
associated with a reduced word $\textbf{i}$ of $w_0$,
we get a free crystal $X_*(T\times (\mathbb{C}^{\times})^{l(w_0)})$ by the tropicalization functor.
Considering the tropicalization of a rational function $\Phi_{BK}$ on $T\cdot B_{w_0}^-$, 
one obtain a subcrystal 
\[
\{ z\in X_*(T\times (\mathbb{C}^{\times})^{l(w_0)}) | {\rm Trop}(\Phi_{BK}\circ \theta_{\textbf{i}})(z)\geq0 \},
\]
which is isomorphic to the disjoint union of all crystal bases $B(\lambda)$ of
the finite dimensional irreducible representations of $U_q(^L\mathfrak{g})$ with highest weights $\lambda$ \cite{BK}. 
Here $^L\mathfrak{g}$ is the Langlands dual Lie algebra of $\mathfrak{g}={\rm Lie}(G)$.
The function $\Phi_{BK}$ is defined as
\begin{equation}\label{BK-def}
\Phi_{BK}=\sum_{i\in I}\frac{\Delta_{w_0\Lambda_i,s_i\Lambda_i}}{\Delta_{w_0\Lambda_i,\Lambda_i}}
+
\sum_{i\in I}\frac{\Delta_{w_0s_i\Lambda_i,\Lambda_i}}{\Delta_{w_0\Lambda_i,\Lambda_i}}
\end{equation}
and called a {\it Berenstein-Kazhdan decoration function} or {\it BK potential function}.
Here, $\Lambda_i$ is the $i$-th fundamental weight, for $u$, $v\in W$, a function $\Delta_{u\Lambda_i,v\Lambda_i}$
is a generalized minor (Definition \ref{gen-def}).

In \cite{KN}, a {\it half potential} $\Phi^h_{BK}=\sum_{i\in I}\Delta_{w_0\Lambda_i,s_i\Lambda_i}$ is introduced,
which is the restriction of the first term of $\Phi_{BK}$ in (\ref{BK-def}) to $B_{w_0}^-$.
Defining the positive structure $\theta^-_{\textbf{i}}: (\mathbb{C}^{\times})^{l(w_0)}\rightarrow B_{w_0}^-$ on $B_{w_0}^-$,
we get a subcrystal
\begin{equation}\label{cone1}
\{ z\in X_*(\mathbb{C}^{\times})^{l(w_0)}) | {\rm Trop}(\Phi^h_{BK}\circ \theta_{\textbf{i}}^-)(z)\geq0 \},
\end{equation}
which is isomorphic to the crystal base $B(\infty)$ of the negative part $U_q^-(^L\mathfrak{g})$.
Finding monomials appearing in $\Phi^h_{BK}\circ \theta_{\textbf{i}}^-$, one obtain an explicit form of $B(\infty)$
in $X_*((\mathbb{C}^{\times})^{l(w_0)})\cong \mathbb{Z}^{l(w_0)}$. 
It is important to find the explicit form
for studies of the string cones or polyhedral realizations of $B(\infty)$ \cite{Lit,NZ}.
As related works,
in the case $\mathfrak{g}$ is of type ${\rm A}_n$, Gleizer and Postnikov
gave a purely combinatorial rule to compute the inequalities which define the string cone
using the rigorous paths in a graph constructed from wiring diagrams \cite{GP}.
In \cite{GKS16}, a combinatorial expression of string cones via dual Reineke vectors constructed by rhombus tiling tools
is given. 

In this paper, firstly, we invent an algorithm for computing 
$\Delta_{w_0\Lambda_i,s_i\Lambda_i}\circ \theta_{\textbf{i}}^-$ explicitly in the case $V(\Lambda_i)$
is a minuscule representation (Theorem \ref{thm1}). 
The algorithm computes a graph $DG$ whose vertices are labelled by all monomials in 
$\Delta_{w_0\Lambda_i,s_i\Lambda_i}\circ \theta_{\textbf{i}}^-(t_1,\cdots,t_{l(w_0)})$.
From information of exponents of a monomial, 
our algorithm generates new monomials one after another
without any complicated combinatorial tools.
The rule to construct the new monomials is similar to the computation of polyhedral realizations given in \cite{NZ}
or action of Kashiwara operators of {\em monomial realizations}
of Kashiwara's crystals \cite{K,Nj}.
In particular, if $\mathfrak{g}$ is of type ${\rm A}_n$ then we calculate
$\Delta_{w_0\Lambda_i,s_i\Lambda_i}$ for all $i\in I$ so that one can compute the explicit form of (\ref{cone1}).
For example, if $\mathfrak{g}$ is of type ${\rm A}_3$ and $\textbf{i}=(1,2,3,2,1,2)$ then our algorithm computes
the following graphs $DG$ of $\Delta_{w_0\Lambda_i,s_i\Lambda_i}$:
\[
DG\ (i=1):\ \ t_4\rightarrow \frac{t_5}{t_6},
\]
\[
DG\ (i=2):\ \ t_6,
\]
\[
DG\ (i=3):\ \ t_1\rightarrow \frac{t_2t_4}{t_5}\rightarrow \frac{t_2}{t_6}\rightarrow \frac{t_3}{t_4t_6}.
\]
By these graphs, it follows
\[
\Delta_{w_0\Lambda_1,s_1\Lambda_1}\circ\theta^-_{\textbf{i}}(t_1,t_2,\cdots,t_6)=t_4+\frac{t_5}{t_6},\quad
\Delta_{w_0\Lambda_2,s_2\Lambda_2}\circ\theta^-_{\textbf{i}}(t_1,t_2,\cdots,t_6)=t_6,
\]
\[
\Delta_{w_0\Lambda_3,s_3\Lambda_3}\circ\theta^-_{\textbf{i}}(t_1,t_2,\cdots,t_6)=
t_1+ \frac{t_2t_4}{t_5}+ \frac{t_2}{t_6}+ \frac{t_3}{t_4t_6}.
\]
Considering the tropicalization functor, one obtain the explicit form of (\ref{cone1}):
\[
B(\infty)
\cong
\left\{ (z_1,\cdots,z_6)\in X_*(\mathbb{C}^{\times})^{6})\cong\mathbb{Z}^6 \left| 
\begin{array}{l}
 z_4\geq0, z_5-z_6\geq0, z_6\geq0, \\
 z_1\geq0,
 z_2+z_4-z_5\geq0,\\
z_2-z_6\geq0, z_3-z_4-z_6\geq0
\end{array}\right.
\right\}.
\]
Next, we prove that our graph $DG$ has some properties similar to crystal graphs of minuscule representations.
The monomials in $DG$ are parametrized by $\textbf{i}$-trails (\cite{BZ}, see Theorem \ref{trail-thm})
and a bijection between a set of $\textbf{i}$-trails of bipartite type and a Demazure crystal of a 
fundamental representation are given in \cite{J}, which gives us a motivation to compare the graph $DG$
and crystal graphs.

Finally, we show that our algorithm works for the case $\mathfrak{g}$ is of type ${\rm G}_2$ though the representations $V(\Lambda_i)$ ($i=1,2$)
are not minuscule. Thus, one can compute the explicit form of (\ref{cone1}).
We also give other examples such that $V(\Lambda_i)$ is non-minuscule and our algorithm works.
From these examples, in the case $\textbf{i}$ is `adapted', that is, it satisfies the assumption in Remark \ref{last-rem},
we expect our algorithm works and the monomials in the graph $DG$ appear in a crystal graph of a monomial realization for a
fundamental representation.

The organization of this paper is as follows:
In Section 2, we review generalized minors and their calculation by using $\textbf{i}$-trails.
Section 3 provides a background on geometric crystals and decoration functions.
In Section 4, several properties of generalized minors $\Delta_{w_0\Lambda_i,s_i\Lambda_i}$ are shown.
Section 5 presents our main results. We give an algorithm to compute the graph $DG$
when $V(\Lambda_i)$ is a minuscule representation in Theorem \ref{thm1}. 
Next, we prove the graph $DG$ has several similar properties to crystal graphs of minuscule representations.
In Section 6, we consider the case $V(\Lambda_i)$ is non-minuscule.
It is shown the algorithm in Theorem \ref{thm1} works in the case $G$ is of type ${\rm G}_2$ via
direct computations. Furthermore, we give several examples the algorithm works when $V(\Lambda_i)$ is non-minuscule
and $\textbf{i}$ is a specific one.

\vspace{2mm}

\nd
{\bf Acknowledgements}
Y.K. is supported by JSPS KAKENHI Grant Number JP20J00186.
G.K. is partially supported by Laboratory of Mirror Symmetry NRU HSE, RF Government grant, ag. N 14.641.31.0001.
T.N. was supported in part by JSPS KAKENHI Grant Number JP20K03564.
Y.K. and T.N. would like to thank to MCCME and V.Poncelet Center in Moscow for hospitality.
G.K. also thanks Sophia University for hospitality.
We thank Denis Mironov for writing a code in Sage based on the algorithm in Theorem \ref{thm1}.

\section{Trails and generalized minors}

\subsection{Notation}\label{notation}

Let $G$ be a simply connected connected simple algebraic group,
$B,\ B^-\subset G$ Borel subgroups, $T:=B\cap B^-$ the maximal torus,
$W={\rm Norm}_G(T)/T$ Weyl group, 
$U$, $U^-$ be unipotent radicals
of $B$, $B^-$,
$A=(a_{i,j})$ the Cartan matrix
of $G$ with an index set $I=\{1,2,\cdots,n\}$. We set $\mathfrak{g}={\rm Lie}(G)$ with Chevalley generators
$e_i$, $f_i$, $h_i$ ($i\in I$), a Cartan subalgebra $\mathfrak{h}$ and the pairing $\langle, \rangle$ of $\mathfrak{h}$ and $\mathfrak{h}^*$.
We use the same numbering of the vertices in the Dynkin diagram as in Kac's book \cite{Kac}.
Let $\Lambda_i$ be the $i$-th fundamental weight, that is, $\langle h_j, \Lambda_i\rangle=\delta_{j,i}$.
Let $P=\oplus_{i\in I}\mathbb{Z}\Lm_i$ be the weight lattice, 
$P_+=\oplus_{i\in I}\mathbb{Z}_{\geq0}\Lm_i$ the positive weight lattice,
 $P^*=\oplus_{i\in I}\mathbb{Z}h_i$ the dual weight lattice,
$\{\alpha_i\}$ ($i\in I$) the set of simple roots.
For $\lambda\in P_+$, let $V(\lambda)$ denote the finite dimensional irreducible $\mathfrak{g}$-module with highest weight $\lambda$.
For two integers $l$, $m\in\mathbb{Z}$ with $l\leq m$, we set $[l,m]:=\{l,l+1,\cdots,m-1,m\}$.

\subsection{An open embedding}\label{aoe}

For $i\in I$ and $t\in\mathbb{C}$, we set
\[
x_i(t):={\rm exp}(te_i),\ y_i(t):={\rm exp}(tf_i)\in G.
\]
Let $\phi_i : SL_{2}(\mathbb{C})\rightarrow G$ be the canonical embedding such that
\[
x_i(t)=\phi_i\left(
\begin{pmatrix}
1 & t \\
0 & 1
\end{pmatrix}
\right),\quad
y_i(t)=\phi_i\left(
\begin{pmatrix}
1 & 0 \\
t & 1
\end{pmatrix}
\right).
\]
For $i\in I$ and $t\in\mathbb{C}^{\times}$,
one sets
\[
t^{h_i}:=
\phi_i\left(
\begin{pmatrix}
t & 0 \\
0 & t^{-1}
\end{pmatrix}
\right)\in T
\]
and
\[
x_{-i}(t):=y_i(t)t^{-h_i}=
\phi_i\left(
\begin{pmatrix}
t^{-1} & 0 \\
1 & t
\end{pmatrix}
\right)\in G.
\]
We also define
\[
\overline{s_i}:=x_i(-1)y_i(1)x_i(-1)\in {\rm Norm}_G(T),
\]
which is a representative of a simple reflection $s_i\in W={\rm Norm}_G(T)/T$.
For each $w\in W$, one can define the a representative $\overline{w}\in{\rm Norm}_G(T)$ by the rule
$\overline{uv}=\overline{u}\cdot \overline{v}$ if $l(uv)=l(u)+l(v)$. We put $B^-_w:=B^-\cap U\overline{w}U $ for $w\in W$.
For a reduced word $\textbf{i}=(i_1,\cdots,i_N)$ of the longest element $w_0\in W$, one defines a map
$\theta^-_{\textbf{i}}: (\mathbb{C}^{\times})^{N}\rightarrow G$
as
\[
\theta^-_{\textbf{i}}(t_1,\cdots,t_N):=x_{-i_1}(t_1)\cdots x_{-i_N}(t_N).
\]
\begin{prop}\label{OEprop}\cite{BZ}
$\theta^-_{\textbf{i}}$ is an open embedding $(\mathbb{C}^{\times})^{N}\hookrightarrow B^-_{w_0}$.
\end{prop}

\subsection{Generalized minors}

Let $G_0:=U^-TU$ be the open subset of $G$ whose elements $x\in G_0$ are uniquely decomposed as
$x=[x]_-[x]_0[x]_+$ with $[x]_-\in U^-$, $[x]_0\in T$ and $[x]_+\in U$. 
\begin{defn}\label{gen-def}\cite{FZ}
For $u,v\in W$ and $i\in I$, the generalized minor $\Delta_{u\Lambda_i,v\Lambda_i}$ is the regular function on $G$
such that for any $x\in\overline{u}G_0\overline{v}^{-1}$, it holds
\[
\Delta_{u\Lambda_i,v\Lambda_i}(x)=([\overline{u}^{-1}x\overline{v}]_0)^{\Lambda_i}.
\]
Here, for $t\in \mathbb{C}^{\times}$ and $j\in I$, we set $(t^{h_j})^{\Lambda_i}=(t^{\Lambda_i(h_j)})$ and extend it to
the group homomorphism $T\rightarrow \mathbb{C}^{\times}$.
\end{defn}

We can compute generalized minors by a representation theoretical way.
First, let $\omega:\mathfrak{g}\to\mathfrak{g}$ be the anti-involution 
\[
\omega(e_i)=f_i,\quad
\omega(f_i)=e_i,\quad \omega(h)=h\ (h\in P^*)
\]
and we extend it to $G$ by setting
$\omega(x_i(c))=y_{i}(c)$, $\omega(y_{i}(c))=x_i(c)$ and $\omega(t)=t$
for $i\in I$ and $t\in T$. Here, $x_i(c)$ and $y_i(c)$ are defined in subsection \ref{aoe}.
Note that $\omega(\ovl s_i^{\pm1})=\ovl s_i^{\mp1}$.
There exists a $\mathfrak{g}$ (or $G$)-invariant bilinear form on the
finite-dimensional irreducible
$\mathfrak{g}$-module $V(\Lm_i)$ with a highest weight vector $u_{\Lambda_i}$ such that
$\langle u_{\Lambda_i},u_{\Lambda_i}\rangle=1$ and
\[
 \langle au,v\rangle=\langle u,\omega(a)v\rangle,
\quad(u,v\in V(\Lm_i),\,\, a\in \mathfrak{g}\ (\text{or }G)).
\]
For $g\in G$, 
we have the following simple fact:
\[
 \Delta_{\Lm_i,\Lm_i}(g)=\langle gu_{\Lm_i},u_{\Lm_i}\rangle.
\]
Thus, for $w,\ w'\in W$, one obtain
\begin{equation}\label{minor-bilin}
 \Delta_{w'\Lm_i,w\Lm_i}(g)=
\Delta_{\Lm_i,\Lm_i}(\ovl{w'}^{-1}g\ovl w)=
\langle g\ovl w\cdot u_{\Lm_i},\ovl{w'}\cdot u_{\Lm_i}\rangle.
\end{equation}

\subsection{\textbf{i}-trails}

\begin{defn}\label{pretrail}
For a finite dimensional representation $V$ of $\mathfrak{g}$, two weights $\gamma$, $\delta$ of $V$ and
a sequence $\textbf{i}=(i_1,\cdots,i_l)$ of indices from $[1,n]$,
we say a sequence
$\pi=(\gamma=\gamma_0,\gamma_1,\cdots,\gamma_l=\delta)$ 
is a {\it pre-\textbf{i}-trail} from $\gamma$ to $\delta$ if
for $k\in[1,l]$, it holds $\gamma_k\in P$ and $\gamma_{k-1}-\gamma_k=c_k\alpha_{i_k}$ with some nonnegative integer $c_k$.
\end{defn}
\nd
Note that for $k\in[1,N]$, it holds
\begin{equation}\label{ck}
c_k=\frac{\gamma_{k-1}-\gamma_k}{2}(h_{i_k}).
\end{equation}

\begin{defn}\label{def-tr}\cite{BZ}
In the setting of Definition \ref{pretrail},
if pre-$\textbf{i}$-trail $\pi$ satisfies
\begin{itemize}
\item $e^{c_1}_{i_1}e^{c_2}_{i_2}\cdots e^{c_l}_{i_l}$ is a non-zero linear map from $V_{\delta}$ to $V_{\gamma}$,
\end{itemize}
then $\pi$ is called an {\it \textbf{i}-trail} from $\gamma$ to $\delta$,
where $V=\oplus_{\mu} V_{\mu}$ is the weight decomposition of $V$.
\end{defn}

\nd
For a pre-$\textbf{i}$-trail $\pi$, we set
\begin{equation}\label{dk}
d_k(\pi):=\frac{\gamma_{k-1}+\gamma_k}{2}(h_{i_k}).
\end{equation}
Note that $d_k(\pi)=c_k+\gamma_k(h_{i_k})\in\mathbb{Z}$ by (\ref{ck}).
If $\gamma_{k-1}=s_{i_k}\gamma_k$ then $d_k(\pi)=0$.
For a reduced word $\textbf{i}=(i_1,\cdots,i_l)$ of an element in $W$, we define
$\theta^{-}_{\textbf{i}}(t_1,\cdots,t_l):=x_{-i_1}(t_1)\cdots x_{-i_l}(t_l)\in G$, which coincides
with the open embedding in Proposition \ref{OEprop} when $w=w_0$.

\begin{thm}\label{trail-thm}\cite{BZ}
For $u$, $v\in W$ and $i\in[1,n]$, we see that
$\Delta_{u\Lambda_i,v\Lambda_i}(\theta^-_{\textbf{i}}(t_1,\cdots,t_l))$
is a linear combination of the monomials $t_1^{d_1(\pi)}\cdots t_l^{d_l(\pi)}$ with positive coefficients for
all $\textbf{i}$-trail from $-u\Lambda_i$ to $-v\Lambda_i$ in $V(-w_0\Lambda_i)$.
\end{thm}

Let $\pi=(\gamma_0,\gamma_1,\cdots,\gamma_l)$ and
$\pi'=(\gamma_0',\gamma_1',\cdots,\gamma_l')$ be two $\textbf{i}$-trails
 from $-u\Lambda_i$ to $-v\Lambda_i$ with integers $\{c_k\}_{k\in[1,l]}$, $\{c_k'\}_{k\in[1,l]}$
 such that $\gamma_{k-1}-\gamma_k=c_k\alpha_{i_k}$, $\gamma_{k-1}'-\gamma_k'=c_k'\alpha_{i_k}$.
If $d_k(\pi)=d_k(\pi')$ for all $k\in[1,l]$ then
by the relation $d_k(\pi)=c_k+\gamma_k(h_{i_k})$ and $\gamma_l=\gamma_l'=-v\Lambda_i$,
one can inductively show $\gamma_k=\gamma_{k}'$ and $c_k=c_{k}'$ for $k=l, l-1,\cdots,1$ so that $\pi=\pi'$.
Therefore, for each monomial $M$ in $\Delta_{u\Lambda_i,v\Lambda_i}(\theta^-_{\textbf{i}}(t_1,\cdots,t_l))$,
there uniquely exists a corresponding $\textbf{i}$-trail $\pi$ such that $M=t_1^{d_1(\pi)}\cdots t_l^{d_l(\pi)}$.

\section{Decorated geometric crystals}

\subsection{Crystals and monomial realizations}\label{cry-mono}

First, let us recall the definition of crystals. We use
a slightly different notation from \cite{K2}.

\begin{defn}\cite{K2}
A $\mathfrak{g}$-{\it crystal} is a set $\mathcal{B}$ together with the maps
$\tilde{\gamma}_i:\mathcal{B}\rightarrow \mathbb{Z}$,
$\tilde{\varepsilon}_i,\tilde{\varphi}_i:\mathcal{B}\rightarrow \mathbb{Z}\sqcup \{-\infty\}$
and $\tilde{e}_i$,$\tilde{f}_i:\mathcal{B}\rightarrow \mathcal{B}\sqcup\{0\}$
($i\in I$) satisfying the followings: For $b,b'\in\mathcal{B}$, $i,j\in I$,
\begin{enumerate}
\item[$(1)$] $\tilde{\varphi}_i(b)=\tilde{\varepsilon}_i(b)+\tilde{\gamma}_i(b)$,
\item[$(2)$] $\tilde{\gamma}_j(\tilde{e}_ib)=\tilde{\gamma}_j(b)+a_{j,i}$ if $\tilde{e}_i(b)\in\mathcal{B}$,
\quad $\tilde{\gamma}_j(\tilde{f}_ib)=\tilde{\gamma}_j(b)-a_{j,i}$ if $\tilde{f}_i(b)\in\mathcal{B}$,
\item[$(3)$] $\tilde{\varepsilon}_i(\tilde{e}_i(b))=\tilde{\varepsilon}_i(b)-1,\ \ 
\tilde{\varphi}_i(\tilde{e}_i(b))=\tilde{\varphi}_i(b)+1$\ if $\tilde{e}_i(b)\in\mathcal{B}$, 
\item[$(4)$] $\tilde{\varepsilon}_i(\tilde{f}_i(b))=\tilde{\varepsilon}_i(b)+1,\ \ 
\tilde{\varphi}_i(\tilde{f}_i(b))=\tilde{\varphi}_i(b)-1$\ if $\tilde{f}_i(b)\in\mathcal{B}$, 
\item[$(5)$] $\tilde{f}_i(b)=b'$ if and only if $b=\tilde{e}_i(b')$,
\item[$(6)$] if $\tilde{\varphi}_i(b)=-\infty$ then $\tilde{e}_i(b)=\tilde{f}_i(b)=0$.
\end{enumerate}
The {\it crystal graph} of $\mathcal{B}$ is the colored directed graph such that vertices are elements in $\mathcal{B}$
and arrows are determined by $b\overset{i}{\rightarrow} b'$ if and only if
$b'=\tilde{f}_ib$ for $b,b'\in\mathcal{B}$, $i\in I$.
We call $\tilde{e}_i$,$\tilde{f}_i$ {\it Kashiwara operators}.
A crystal $\mathcal{B}$ is said to be {\it free} if the Kashiwara operators
$\tilde{e}_i$ $(i\in I)$ are bijections $\tilde{e}_i:\mathcal{B}\rightarrow \mathcal{B}$. 
\end{defn}

It is well-known that the irreducible integrable highest weight module $V_q(\lambda)$ $(\lambda\in P_+)$
of quantized universal enveloping algebra
$U_q(\mathfrak{g})$ has the crystal base $(L(\lambda),B(\lambda))$ and the negative part $U_q^-(\mathfrak{g})$ of $U_q(\mathfrak{g})$ also 
has the crystal base $(L(\infty),B(\infty))$. The sets $B(\lambda)$ and $B(\infty)$ are important and interesting crystals.

Here, we review the {\it monomial realization} of $B(\lambda)$ \cite{K,Nj}.
We fix $p$ a set of integers $\{p_{j,k} | j,k \in I, a_{j,k}<0\}$ such that $p_{j,k}+p_{k,j}=1$.
For doubly-indexed variables $\{Y_{s,i} \,|\, i \in I$, $s\in \mathbb{Z}\}$, let us define a crystal structure on the set of Laurent monomials
\begin{gather*}
{\mathcal Y}:=\left\{Y=\prod\limits_{s \in \mathbb{Z},\ i \in I}
Y_{s,i}^{\zeta_{s,i}}\, \Bigg| \,\zeta_{s,i} \in \mathbb{Z},\
~\text{only finitely many } \zeta_{s,i} \neq 0 \right\}.
\end{gather*}
For $Y=\prod\limits_{s \in \mathbb{Z},\; i \in I} Y_{s,i}^{\zeta_{s,i}}\in {\mathcal Y}$, we set
$\tilde{\gamma}_i(Y):= \sum_{s\in\mathbb{Z}}\zeta_{s,i}$ and
\begin{gather*}
\tilde{\varphi}_i(Y):=\max\left\{\! \sum\limits_{k\leq s}\zeta_{k,i} \,|\, s\in \mathbb{Z} \!\right\},\!
\quad
\tilde{\varepsilon}_i(Y):=\tilde{\varphi}_i(Y)-\tilde{\gamma}_i(Y)
\end{gather*}
for $i\in I$. Putting
\[
A_{s,k}:=Y_{s,k}Y_{s+1,k}\prod\limits_{j; a_{j,k}<0}Y_{s+p_{j,k},j}^{a_{j,k}} \quad \text{for }s\in\mathbb{Z},\ k\in I,
\]
one define actions of Kashiwara operators as follows:
\begin{gather*}
\tilde{f}_iY:=
\begin{cases}
A_{n_{f_i},i}^{-1}Y & \text{if} \quad \tilde{\varphi}_i(Y)>0,
\\
0 & \text{if} \quad \tilde{\varphi}_i(Y)=0,
\end{cases}
\quad
\tilde{e}_iY:=
\begin{cases}
A_{n_{e_i},i}Y & \text{if} \quad \tilde{\varepsilon}_i(Y)>0,
\\
0 & \text{if} \quad \tilde{\varepsilon}_i(Y)=0,
\end{cases}
\end{gather*}
where
\begin{gather*}
n_{f_i}:=\min \left\{r\in\mathbb{Z} \,\Bigg|\, \tilde{\varphi}_i(Y)= \sum\limits_{k\leq r}\zeta_{k,i}\right\},
\qquad
n_{e_i}:=\max \left\{r\in\mathbb{Z} \,\Bigg|\, \tilde{\varphi}_i(Y)= \sum\limits_{k\leq r}\zeta_{k,i}\right\}.
\end{gather*}

\begin{thm}\cite{K,Nj}
\begin{enumerate}
\item[(i)] The set ${\mathcal Y}$ together with the above maps
$\tilde{\gamma}_i$, $\tilde{\varepsilon}_i$,  $\tilde{\varphi}_i$, and $\tilde{e}_i$, $\tilde{f}_i$ $(i\in I)$ is a crystal.
\item[(ii)] If a~monomial $Y \in {\mathcal Y}$ satisfies $\tilde{\varepsilon}_i(Y)=0$ for all $i \in I$ then the set
\[
\{\tilde{f}_{j_s}\cdots\tilde{f}_{j_1}Y | s\in\mathbb{Z}_{\geq0},\ j_1,\cdots,j_s \in I \}\setminus\{0\}
\]
is isomorphic to the crystal $B(\sum_{i\in I}\tilde{\gamma}_i(Y)\Lambda_i)$.
\end{enumerate}
\end{thm}

\subsection{Geometric crystals and tropicalization functor}\label{GC-trop}

Next, we review the definition of geometric crystals, which is a
geometric analog of crystals.

\begin{defn}\cite{BK}
Let $X$ be an irreducible algebraic variety over $\mathbb{C}$. We suppose that
for each $i\in I$, rational functions $\gamma_i$, $\varepsilon_i$ on $X$
and a unital rational $\mathbb{C}^{\times}$ action
$\ovl{e}_i : \mathbb{C}^{\times}\times X\rightarrow X$ on $X$ are defined.
The quadruple $(X,\{\ovl{e}_i\}_{i\in I}, \{\gamma_i\}_{i\in I}, \{\varepsilon_i\}_{i\in I})$ is said to be a
$\mathfrak{g}$-{\it geometric crystal} if the following holds:
\begin{enumerate}
\item $(\{1\}\times X)\cap {\rm dom}(\ovl{e}_i)$ is open dense in $\{1\}\times X$ for each $i\in I$, where
${\rm dom}(\ovl{e}_i)$
denotes the domain of definition of $\ovl{e}_i : \mathbb{C}^{\times}\times X\rightarrow X$.
We write $\ovl{e}_i^c(x):=\ovl{e}_i(c,x)$ for $c\in \mathbb{C}^{\times}$ and $x\in X$.
\item It holds $\gamma_j(\ovl{e}_i^c(x))=c^{a_{i,j}}\gamma_j(x)$ for any $i,j\in I$ and $c\in \mathbb{C}^{\times}$.
\item For $i,j\in I$ and $c_1,c_2\in\mathbb{C}^{\times}$,
\begin{eqnarray*}
\ovl{e}_i^{c_1}\ovl{e}_j^{c_2}&=&\ovl{e}_j^{c_2}\ovl{e}_i^{c_1} \quad \text{ if }a_{i,j}=0,\\
\ovl{e}_i^{c_1}\ovl{e}_j^{c_1c_2}\ovl{e}_i^{c_2}
&=&\ovl{e}_j^{c_2}\ovl{e}_i^{c_1c_2}\ovl{e}_j^{c_1} \quad \text{ if }a_{i,j}=a_{j,i}=-1,\\
\ovl{e}_i^{c_1}\ovl{e}_j^{c_1^2c_2}\ovl{e}_i^{c_1c_2}\ovl{e}_j^{c_2}
&=&\ovl{e}_j^{c_2}\ovl{e}_i^{c_1c_2}\ovl{e}_j^{c_1^2c_2}\ovl{e}_i^{c_1} \quad \text{ if }a_{i,j}=-2,\ a_{j,i}=-1,\\
\ovl{e}_i^{c_1}\ovl{e}_j^{c_1^3c_2}\ovl{e}_i^{c_1^2c_2}\ovl{e}_j^{c_1^3c_2^2}\ovl{e}_i^{c_1c_2}\ovl{e}_j^{c_2}
&=&\ovl{e}_j^{c_2}\ovl{e}_i^{c_1c_2}\ovl{e}_j^{c_1^3c_2^2}\ovl{e}_i^{c_1^2c_2}\ovl{e}_j^{c_1^3c_2}\ovl{e}_i^{c_1} \quad 
\text{ if }a_{i,j}=-3,\ a_{j,i}=-1.
\end{eqnarray*}
\item For any $i,j\in I$ and $c\in\mathbb{C}^{\times}$,
\[
\varepsilon_i(\ovl{e}_i^c(x))=c^{-1}\varepsilon_i(x),\ \quad
\varepsilon_i(\ovl{e}_j^c(x))=\varepsilon_i(x)\ \text{ if }a_{i,j}=0.
\]
\end{enumerate}
We also set $\varphi_i(x):=\gamma_i(x)\varepsilon_i(x)$ for $i\in I$.
\end{defn}

\begin{rem}\label{fliprem}
The above definitions of $\varepsilon_i$, $\varphi_i$ are different from \cite{BK}.
By the flip $\varepsilon_i\mapsto \varepsilon_i^{-1}$, $\varphi_i\mapsto \varphi_i^{-1}$, they coincide
with the above one. We follow the notation in \cite{KN}.
\end{rem}

\begin{defn}
Let $\chi=(X,\{\ovl{e}_i\}_{i\in I}, \{\gamma_i\}_{i\in I}, \{\varepsilon_i\}_{i\in I})$ be
a $\mathfrak{g}$-geometric crystal. For a rational function $f:X\rightarrow \mathbb{C}$,
the pair $(\chi,f)$ is said to be a $\mathfrak{g}$-upper (resp. lower) half-decorated geometric crystal if
\[
f(\ovl{e}_i^c(x))=f(x)+(c^{-1}-1)\varepsilon_i(x)\quad
(\text{resp. } f(\ovl{e}_i^c(x))=f(x)+(c-1)\varphi_i(x))
\]
for any $i\in I$, $c\in\mathbb{C}^{\times}$ and $x\in X$. Then $f$ is said to be an
upper (resp. lower) half-decoration or upper (resp. lower) half-potential.
\end{defn}

Let $^t\mathbb{Z}:=\mathbb{Z}\sqcup\{-\infty\}$ be the tropical semi-field with the multiplication `$+$'
and summation `min'. We define a map $\mathbb{V}:\mathbb{C}(x)\rightarrow$$^t\mathbb{Z}$ as
\[
\mathbb{V}(f(x)):=
\begin{cases}
-{\rm deg}(f(x^{-1})) & {\rm if}\ f\not\equiv0,\\
-\infty & {\rm if}\ f\equiv0,
\end{cases}
\]
where for $g(x)=\sum^m_{j=0}a_jx^j\in\mathbb{C}[x]$ and $h(x)=\sum^{m'}_{j=0}b_jx^j\in\mathbb{C}[x]\setminus \{0\}$
with $a_m\neq0$, $b_m\neq0$,
we set ${\rm deg}(g(x)/h(x)):=m-m'$. For an algebraic torus $T'$ over $\mathbb{C}$, let $X^*(T'):={\rm Hom}(T',\mathbb{C}^{\times})$ and
$X_*(T'):={\rm Hom}(\mathbb{C}^{\times},T')$ be the sets of characters and cocharacters, respectively. 
\begin{defn}
\begin{enumerate}
\item[$(i)$] Let $T'$ be an algebraic torus. 
A rational function $f:T'\rightarrow \mathbb{C}$ is said to be positive
if $f$ can be written as
\[
f=\frac{g}{h}
\]
with some regular functions $g=\sum_{\mu\in X^*(T')}a_{\mu}\mu$ and 
$h=\sum_{\mu\in X^*(T')}b_{\mu}\mu (\neq0)$
such that all $a_{\mu},\ b_{\mu}$ are nonnegative integers. 
\item[$(ii)$]
Let $f:T'\rightarrow T''$ be a rational map between
two algebraic tori $T'$, $T''$. The map $f$ is said to be positive
if $\xi\circ f : T'\rightarrow \mathbb{C}$ is a positive rational function for any $\xi\in X^*(T'')$. 
\end{enumerate}
\end{defn}
For two positive rational functions $f_1$ and $f_2\in\mathbb{C}(x)\setminus\{0\}$ on $\mathbb{C}^{\times}$, it holds
\[
\mathbb{V}(f_1\cdot f_2)=\mathbb{V}(f_1)+\mathbb{V}(f_2),
\qquad
\mathbb{V}(f_1/ f_2)=\mathbb{V}(f_1)-\mathbb{V}(f_2),
\]
\[
\mathbb{V}(f_1+ f_2)={\rm min}(\mathbb{V}(f_1),\mathbb{V}(f_2)).
\]
\begin{defn}
Let $T'$, $T''$ be algebraic tori, $f:T'\rightarrow T''$ be a rational map. 
Let $\langle,\rangle$ be the pairing between $X^*(T'')$ and $X_*(T'')$.
We define a map
$\widehat{f}:X_*(T')\rightarrow X_*(T'')$ as for $\chi\in X^*(T'')$ and $\xi\in X_*(T')$,
\[
\langle\chi,\widehat{f}(\xi)\rangle = \mathbb{V}(\chi\circ f\circ\xi).
\]
\end{defn} 
Let $\mathcal{T}_+$ be the category whose objects are algebraic tori over $\mathbb{C}$
and morphisms are positive rational maps and
$\mathfrak{S}\mathfrak{e}\mathfrak{t}$ be the category of sets.
We can define a functor ${\rm Trop} : \mathcal{T}_+ \rightarrow \mathfrak{S}\mathfrak{e}\mathfrak{t}$
as
\[
T'\mapsto X_*(T'), \quad (f:T'\rightarrow T'') \mapsto (\widehat{f} : X_*(T')\rightarrow X_*(T'')).
\]
The functor ${\rm Trop}$ is called a {\it tropicalization}. For instance, if $T'=(\mathbb{C}^{\times})^3$,
$T''=\mathbb{C}^{\times}$ and $f:T'\rightarrow T''$ is defined as
\[
f(x_1,x_2,x_3)=\frac{x_1^2+2x_2x_3+x_3^2}{x_1+x_2}
\]
then ${\rm Trop}(f):X_*((\mathbb{C}^{\times})^3)\rightarrow X_*(\mathbb{C}^{\times})$ is given by
\[
{\rm Trop}(f)(z_1,z_2,z_3)={\rm min}(2z_1,z_2+z_3,2z_3)-{\rm min}(z_1,z_2),
\]
where we identify $X_*((\mathbb{C}^{\times})^3)$, $X_*(\mathbb{C}^{\times})$ with $\mathbb{Z}^3$, $\mathbb{Z}$.
In this way, the product $\times$, sum $+$ and division $\div$ in $f$ correspond to
the sum $+$, ${\rm min}$ and minus $-$ in ${\rm Trop}(f)$ respectively.

\subsection{Tropicalizations of geometric crystals}\label{Trop-geom}

\begin{defn}
Let $\chi=(X,\{\ovl{e}_i\}_{i\in I}, \{\gamma_i\}_{i\in I}, \{\varepsilon_i\}_{i\in I})$ be
a $\mathfrak{g}$-geometric crystal, $T'$ an algebraic torus and $\theta:T'\rightarrow X$ a birational map.
\begin{enumerate}
\item[(i)] For $i\in I$, rational functions $\gamma_i\circ \theta:T'\rightarrow \mathbb{C}$ and
$\varepsilon_i\circ \theta:T'\rightarrow \mathbb{C}$ are positive.
\item[(ii)] For $i\in I$, the rational map $\ovl{e}_{i,\theta}:\mathbb{C}^{\times}\times T'\rightarrow \mathbb{C}^{\times}$
defined as $(c,t)\mapsto \theta^{-1}\circ \ovl{e}_i^c\circ\theta(t)$ is positive.
\end{enumerate}
Then $\theta$ is said to be a {\it positive structure} on $\chi$. Additionally,
if $f$ is an upper or lower half decoration of $\chi$ and the conditions (i), (ii) and 
\begin{enumerate}
\item[(iii)] the rational function $f\circ \theta : T'\rightarrow \mathbb{C}$ is positive
\end{enumerate}
hold then $\theta$ is said to be a positive structure on $(\chi,f)$.
\end{defn}

\begin{thm}\cite{BK0,N0}
Let $\chi=(X,\{\ovl{e}_i\}_{i\in I}, \{\gamma_i\}_{i\in I}, \{\varepsilon_i\}_{i\in I})$ be
a $\mathfrak{g}$-geometric crystal, $\theta:T'\rightarrow X$ be its positive structure.
We set
\[
\tilde{e}_i:={\rm Trop}(\ovl{e}_{i,\theta}):\mathbb{Z}\times X_*(T')\rightarrow X_*(T'),
\]
\[
\tilde{\gamma}_i:={\rm Trop}(\gamma_i\circ \theta):X_*(T')\rightarrow\mathbb{Z},\quad 
\tilde{\varepsilon}_i:={\rm Trop}(\varepsilon_i\circ \theta):X_*(T')\rightarrow\mathbb{Z}.
\]
Then the $6$-tuple $(X_*(T'),\{\tilde{e}_i(1,\cdot)\}_{i\in I},\{\tilde{e}_i(-1,\cdot)\}_{i\in I},\{\tilde{\varepsilon}_i\}_{i\in I},
\{\tilde{\varphi}_i\}_{i\in I},\{\tilde{\gamma}_i\}_{i\in I})$ is a free $^L\mathfrak{g}$-crystal.
Here, $^L\mathfrak{g}$ is the finite dimensional simple Lie algebra whose Cartan matrix
is $^tA$ and $\tilde{\varphi}_i:=\tilde{\gamma}_i+\tilde{\varepsilon}_i$.
\end{thm} 

Let $(\chi,f)$ be
a $\mathfrak{g}$-upper (resp. lower) decorated geometric crystal with a positive
structure $\theta : T'\rightarrow X$. Let us consider a subset
\[
\tilde{B}_{\theta,f}:=\{x\in X_*(T') | {\rm Trop}(f\circ \theta)(x)\geq0\} \subset X_*(T').
\]
Defining $\tilde{e}_i(z,x)=0$ if $\tilde{e}_i(z,x)\notin \tilde{B}_{\theta,f}$ for $z=1,-1$, one can define
a crystal structure on $\tilde{B}_{\theta,f}$ as
\begin{equation}\label{gccrystal}
\mathbb{B}_{\theta,f}:=(\tilde{B}_{\theta,f},\{\tilde{e}_i(1,\cdot)|_{\tilde{B}_{\theta,f}}\}_{i\in I},
\{\tilde{e}_i(-1,\cdot)|_{\tilde{B}_{\theta,f}}\}_{i\in I},\{\tilde{\varepsilon}_i|_{\tilde{B}_{\theta,f}}\}_{i\in I},
\{\tilde{\varphi}_i|_{\tilde{B}_{\theta,f}}\}_{i\in I},
\{\tilde{\gamma}_i|_{\tilde{B}_{\theta,f}}\}_{i\in I}).
\end{equation}

\begin{prop}\cite{KN}
If $f$ is upper (resp. lower) half decoration then $\mathbb{B}_{\theta,f}$ is an upper (resp. lower) normal crystal,
that is,
\[
\tilde{\varepsilon}_i(x)={\rm max}\{n\geq0|\tilde{e}_i(n,x)\neq0\},\quad
(\text{resp. } \tilde{\varphi}_i(x)={\rm max}\{n\geq0|\tilde{e}_i(-n,x)\neq0\}).
\]

\end{prop}

\subsection{Geometric crystal structure on $B^-_{w_0}$}

We use the notation in the subsection \ref{notation}. 
One can define an upper half geometric crystal structure on $B^-_{w_0}=B^-\cap U\overline{w_0}U$ following subsection 2.2 of \cite{BK}. 
For $i\in I$, we define a rational function $\gamma_i:B^-_{w_0}\rightarrow \mathbb{C}^{\times}$ as
\[
\gamma_i: B^-_{w_0}\hookrightarrow B^- \overset{\sim}{\longrightarrow} T\times U
\overset{\text{proj}}{\longrightarrow} T \overset{\alpha_i}{\longrightarrow} \mathbb{C}^{\times}.
\]
In the beginning of subsection \ref{aoe}, we defined the canonical embedding $\phi_i:SL_2(\mathbb{C})\rightarrow G$ for
each $i\in I$. There exists the natural projection $pr_i:B^-\rightarrow B^-\cap \phi_i(SL_2(\mathbb{C}))$. Hence, for any
$x\in B^-_{w_0}$, there uniquely exists 
$v=
\begin{pmatrix}
b_{11} & 0 \\
b_{21} & b_{22}
\end{pmatrix}\in SL_2(\mathbb{C})$
such that
$pr_i(x)=\phi_i(v)$. Now, we define the rational function $\varepsilon_i$ on $B^-_{w_0}$ as
\[
\varepsilon_i(x)=\frac{b_{22}}{b_{21}}.
\]
The above definition is different from \cite{BK} by the reason of flip in Remark \ref{fliprem}.
The unital rational $\mathbb{C}^{\times}$-action $\ovl{e}_i:\mathbb{C}^{\times}\times B^-_{w_0}\rightarrow B^-_{w_0}$
is defined as
\[
\ovl{e}_i(c,x):=
x_i((c-1)\varphi_i(x))\cdot x\cdot x_i((c^{-1}-1)\varepsilon_i(x))
\]
if $\varepsilon_i(x)$ is well-defined, that is, $b_{21}\neq0$ and 
$\ovl{e}_i(c,x):=x$ if $\varepsilon_i(x)$ is not well-defined. Here, $x_i$ is the notation defined in the beginning of subsection \ref{aoe}. 
\begin{prop}\cite{BK}
The quadruple $(B^-_{w_0},\{\ovl{e}_i\}_{i\in I},\{\gamma_i\}_{i\in I},\{\varepsilon_i\}_{i\in I})$
is a $\mathfrak{g}$-geometric crystal.
\end{prop}

We define a regular function on $B^-_{w_0}$ as
\begin{equation}\label{def-hBK}
\Phi^{h}_{\rm BK}=
\sum_{i\in I} \Delta_{w_0\Lm_i,s_i\Lm_i}.
\end{equation}

\begin{prop}\cite{KN}
The function $\Phi^{h}_{\rm BK}$ is an upper half decoration on the geometric crystal $B^-_{w_0}$.
\end{prop}

Recall that we defined an open embedding
$\theta^-_{\textbf{i}}:(\mathbb{C}^{\times})^{N}\hookrightarrow B^-_{w_0}$ in
Proposition \ref{OEprop}, which gives a positive structure of $(B^-_{w_0},\Phi^{h}_{\rm BK})$.
Hence, we get a crystal $\mathbb{B}_{\theta^-_{\textbf{i}},\Phi^{h}_{\rm BK}}$ 
as in (\ref{gccrystal}):
\[
\tilde{B}_{\theta^-_{\textbf{i}},\Phi^{h}_{\rm BK}}:=
\{z\in X_{*}((\mathbb{C}^{\times})^{N})| {\rm Trop}(\Phi^{h}_{\rm BK}\circ \theta^-_{\textbf{i}})(z)\geq0 \},
\]
\[
\mathbb{B}_{\theta^-_{\textbf{i}},\Phi^{h}_{\rm BK}}=
(\tilde{B}_{\theta^-_{\textbf{i}},\Phi^{h}_{\rm BK}},\{\tilde{e}_i(1,\cdot)\}_{i\in I},
\{\tilde{e}_i(-1,\cdot)\}_{i\in I},
\{\tilde{\varepsilon}_i\}_{i\in I},
\{\tilde{\varphi}_i\}_{i\in I},
\{\tilde{\gamma}_i\}_{i\in I}).
\]
Here, we omitted the notation of restrictions $|_{\tilde{B}_{\theta^-_{\textbf{i}},\Phi^{h}_{\rm BK}}}$
for $\tilde{e}_i(1,\cdot)$, $\tilde{e}_i(-1,\cdot)$, $\tilde{\varepsilon}_i$, $\tilde{\varphi}_i$.

\begin{thm}\cite{KN}\label{thm1a}
$\mathbb{B}_{\theta^-_{\textbf{i}},\Phi^{h}_{\rm BK}}$
is a $^L\mathfrak{g}$-crystal isomorphic to the crystal
$B(\infty)$.
\end{thm}

\section{Properties of $\Delta_{w_0\Lm_i,s_i\Lm_i}$}

In this section, we prove some properties of $\Delta_{w_0\Lm_i,s_i\Lm_i}$, which is a summand of $\Phi^{h}_{\rm BK}$.
We fix a reduced word $\textbf{i}=(i_1,i_2,\cdots,i_N)$ of the longest element $w_0\in W$.

\begin{prop}\label{low-prop}\cite{KN}
For $i\in I$, the Laurent polynomial
$\Delta_{w_0\Lm_i,s_i\Lm_i}\circ \theta^-_{\textbf{i}}(t_1,\cdots,t_N)$
has a term
\begin{equation}\label{lowest-term}
t_{J} t_{J+1}^{a_{i_{J+1},i}}\cdots t_{N}^{a_{i_N,i}},
\end{equation}
where $J:={\rm max}\{1\leq k\leq N | i_k=i\}$.
\end{prop}

\begin{prop}\label{high-prop}
For $i\in I$, we suppose that
$s_{i_N}s_{i_{N-1}}\cdots s_{i_{k+1}}\alpha_{i_k}=\alpha_i$ with some $k\in[1,N]$.
Then
the Laurent polynomial
$\Delta_{w_0\Lm_i,s_i\Lm_i}\circ \theta^-_{\textbf{i}}(t_1,\cdots,t_N)$
has a term
\begin{equation}\label{highest-term}
t_k.
\end{equation}
\end{prop}

\nd
{\it Proof.}

Since $s_{i_{k+1}}\cdots s_{i_{N}}\alpha_i=\alpha_{i_k}$, the expression $s_{i_{k+1}}\cdots s_{i_{N}}s_{i}$ is reduced (Lemma 3.11 of \cite{Kac}).
Defining $c_l$ ($k+1\leq l \leq N$) as $c_l={\rm max}\{t\in\mathbb{Z}_{\geq0} | e^{t}_{i_l} e^{c_{l+1}}_{i_{l+1}}\cdots
e^{c_N}_{i_N}v_{-s_i\Lambda_i} \neq0 \}$, we get
\begin{equation}\label{pr3-3-1}
e^{c_{k+1}}_{i_{k+1}}\cdots e^{c_N}_{i_N}v_{-s_i\Lambda_i}=v_{-s_{i_{k+1}\cdots s_{i_N}s_i\Lambda_i}},
\end{equation}
where $v_{-w\Lambda_i}$ denotes a non-zero vector in extremal weight  space $V(-w_0\Lambda_i)_{-w\Lambda_i}$ for $w\in W$.

We also obtain
\begin{eqnarray}
\langle h_{i_k}, -s_{i_{k+1}}\cdots s_{i_N}s_i\Lambda_i \rangle 
&=&-\langle h_{i_k}, s_{i_{k+1}}\cdots s_{i_N}(\Lambda_i-\alpha_i) \rangle \nonumber \\
&=&-\langle h_{i_k}, s_{i_{k+1}}\cdots s_{i_N}\Lambda_i \rangle + \langle h_{i_k}, \alpha_{i_k} \rangle \nonumber\\
&=&-\langle h_{i_k}, s_{i_{k+1}}\cdots s_{i_N}\Lambda_i \rangle + 2 \nonumber\\
&=&-\langle s_{i_{N}}\cdots s_{i_{k+1}}h_{i_k}, \Lambda_i \rangle + 2 \nonumber\\
&=&-\langle h_{i}, \Lambda_i \rangle + 2 \nonumber\\
&=& 1. \label{pr3-3-2}
\end{eqnarray}
By
\[
s_{i_{k+1}}\cdots s_{i_N}s_i\Lambda_i
=s_{i_{k+1}}\cdots s_{i_N}(\Lambda_i-\alpha_i)
=s_{i_{k+1}}\cdots s_{i_N}\Lambda_i - \alpha_{i_k}
\]
and
\[
\langle
h_{i_k}, s_{i_{k+1}}\cdots s_{i_N}\Lambda_i \rangle
=\langle s_{i_{N}}\cdots s_{i_{k+1}}h_{i_k}, \Lambda_i \rangle
=\langle h_i, \Lambda_i \rangle=1,
\]
it holds
\begin{equation}\label{pr3-3-3}
s_{i_{k+1}}\cdots s_{i_N}s_i\Lambda_i=s_{i_k}s_{i_{k+1}}\cdots s_{i_N}\Lambda_i,
\end{equation}
which yields
\begin{equation}\label{pr3-3-4}
v_{-s_{i_{k+1}}\cdots s_{i_N}s_i\Lambda_i}=v_{-s_{i_k}s_{i_{k+1}}\cdots s_{i_N}\Lambda_i}.
\end{equation}
Putting
$c_l$ ($1\leq l \leq k-1$) as $c_l={\rm max}\{t\in\mathbb{Z}_{\geq0} | e^{t}_{i_l} e^{c_{l+1}}_{i_{l+1}}\cdots
e^{c_{k-1}}_{i_{k-1}}v_{-s_{i_k}s_{i_{k+1}}\cdots s_{i_N}\Lambda_i} \neq0 \}$, we get
\begin{equation}\label{pr3-3-5}
e^{c_{1}}_{i_{1}}\cdots e^{c_{k-1}}_{i_{k-1}}v_{-s_{i_k}s_{i_{k+1}}\cdots s_{i_N}\Lambda_i}= v_{-w_0\Lambda_i}.
\end{equation}
Combining (\ref{pr3-3-1}), (\ref{pr3-3-4}) and (\ref{pr3-3-5}) and putting $c_k=0$,
we see that the following sequence $\pi=(\gamma_0,\gamma_1,\cdots,\gamma_N)$ of weights is an \textbf{i}-trail
from $-w_0\Lambda_i$ to $-s_i\Lambda_i$ : 
\[
\gamma_N=-s_i\Lambda_i,\ \gamma_{l-1}=\gamma_l+c_l\alpha_{i_l}\ (1\leq l\leq N).
\]
Since for $l\in[1,N]\setminus\{k\}$ it holds $\gamma_{l-1}=s_{i_{l}}\gamma_{l}$, we obtain
$d_l(\pi)=\frac{\gamma_{l-1}+\gamma_{l}}{2}(h_{i_l})=0$. 
The definition $c_k=0$ means $\gamma_{k-1}=\gamma_k=-s_{i_{k+1}}\cdots s_{i_N}s_i\Lambda_i$ and
it follows from (\ref{pr3-3-2}) that
$d_k(\pi)=\frac{\gamma_{k-1}+\gamma_{k}}{2}(h_{i_k})=1$.

Thus $\Delta_{w_0\Lm_i,s_i\Lm_i}\circ \theta^-_{\textbf{i}}(t_1,\cdots,t_N)$
has a term $t_1^{d_1(\pi)}t_2^{d_2(\pi)}\cdots t_N^{d_N(\pi)}=t_k$ by 
Theorem \ref{trail-thm}. \qed

\section{$\Delta_{w_0\Lm_i,s_i\Lm_i}\circ \theta^-_{\textbf{i}}(t_1,\cdots,t_N)$ for minuscule representations $V(\Lambda_i)$}\label{minus}

We use the same notation as in the previous section.
In this section, we take $i\in I$ such that $V(\Lambda_i)$ is a minuscule representation,
that is, the Weyl group $W$ acts transitively on the weights of $V(\Lambda_i)$.
We give an algorithm to compute $\Delta_{w_0\Lm_i,s_i\Lm_i}\circ \theta^-_{\textbf{i}}(t_1,\cdots,t_N)$ explicitly. 
The table of $i\in I$ is as follows \cite{H}. Note that the numbering of Dynkin diagram is the same one as in \cite{Kac}
as said before:

\begin{table}[h]
  \begin{tabular}{|c|c|c|c|c|c|c|} \hline
  Type of $\mathfrak{g}$ & ${\rm A}_n$ & ${\rm B}_n$ & ${\rm C}_n$ & ${\rm D}_n$ & ${\rm E}_6$ & ${\rm E}_7$ \\ \hline
  $i$ & $1,2,\cdots,n$ & $n$ & $1$ & $1$, $n-1$, $n$ & $1$, $5$ & $6$ \\ \hline
  \end{tabular}
\end{table}

\subsection{Main theorem}

By a property of minuscule representations,
for any weight $\mu$ of $V(\Lambda_i)$ (or $V(-w_0\Lambda_i)$) and $t\in I$, it holds
\begin{equation}\label{mini-pro}
\langle
h_t,\mu
\rangle\in\{1,0,-1\}.
\end{equation}
\nd
For $j\in[1,N]$ such that $j^+\leq N$, we define
\begin{equation}\label{ajdef}
A_{j}:=t_jt_{j^+}\prod_{j<l<j^+}t_l^{a_{i_l,i_j}},
\end{equation}
where $j^+:={\rm min}\{ l\in[1,N] | i_l=i_j,\ l>j \}\cup\{N+1\}$. We also set
$j^-:={\rm max}\{ l\in[1,N] | i_l=i_j,\ l<j \}\cup\{0\}$.
The following theorem is our main theorem:
\begin{thm}\label{thm1}
For $i\in I$, let $k\in[1,N]$ such that
$s_{i_N}s_{i_{N-1}}\cdots s_{i_{k+1}}\alpha_{i_k}=\alpha_i$.
Then the set of monomials appearing in
$\Delta_{w_0\Lm_i,s_i\Lm_i}\circ \theta^-_{\textbf{i}}(t_1,\cdots,t_N)$
coincides with the set of Laurent monomials which label the vertices in the directed graph $DG$ obtained
by the following algorithm:
\begin{enumerate}
\item[(1)] Let $DG_0$ be the graph which has only one vertex $t_k$ and no arrow.
\item[(2)] For each sink $M$ of $DG_l$ and for each $j\in[1,N]$ such that $j^+\leq N$, we add
vertices $M\cdot A_{j}^{-1}$ and arrows $M\rightarrow M\cdot A_{j}^{-1}$ to $DG_l$
if and only if
$M$ has a factor $t_j^{+1}$ and
does not have factors $t_{j^+}^{+1}$.
Let $DG_{l+1}$ be the directed graph obtained from $DG_l$ by this step.
\item[(3)] If a graph $DG_r$ has only one sink 
\begin{equation*}
t_{J} t_{J+1}^{a_{i_{J+1},i}}\cdots t_{N}^{a_{i_N,i}},
\end{equation*}
in (\ref{lowest-term}) of Proposition \ref{low-prop} then
we stop this algorithm and put $DG=DG_r$.
\end{enumerate}
Here, `$M$ has a factor $t_j^{+1}$' means the exponent of $t_j$ is positive in $M$.
\end{thm}
We call the graph $DG$ {\it decoration graph}.
Note that two same monomials obtained from different sinks define the same vertex in step (2).

\vspace{2mm}

\nd
As for the coefficients, we can prove the following:
\begin{prop}\label{coefprop}
The coefficient of each monomial in $\Delta_{w_0\Lm_i,s_i\Lm_i}\circ \theta^-_{\textbf{i}}(t_1,\cdots,t_N)$
is $1$.
\end{prop}

By Theorem \ref{thm1} and Proposition \ref{coefprop}, the sum of all Laurent monomials in the decoration graph $DG$ coincides
with $\Delta_{w_0\Lm_i,s_i\Lm_i}\circ \theta^-_{\textbf{i}}(t_1,\cdots,t_N)$.

\vspace{3mm}

\nd
{\it Proof.}

Using the bilinear form in (\ref{minor-bilin}), it follows
\begin{eqnarray*}
\Delta_{w_0\Lm_i,s_i\Lm_i}\circ \theta^-_{\textbf{i}}(t_1,\cdots,t_N)
&=&
\langle \theta^-_{\textbf{i}}(t_1,\cdots,t_N)
\ovl{s_i}\cdot u_{\Lm_i},\ovl{w_0}\cdot u_{\Lm_i}\rangle\\
&=&\langle
{\rm exp}(t_1f_{i_1})t_1^{-h_{i_1}}\cdots {\rm exp}(t_Nf_{i_N})t_N^{-h_{i_N}}
u_{s_i\Lm_i},u_{w_0\Lm_i}\rangle.
\end{eqnarray*}
Here, $u_{w\Lambda_i}$ denotes the extremal weight vector of weight $w\Lambda_i$ for $w\in W$.
It follows by (\ref{mini-pro}) that $f_j^2=0$ on $V(\Lambda_i)$ for any $j\in I$. Hence,
\[
\Delta_{w_0\Lm_i,s_i\Lm_i}\circ \theta^-_{\textbf{i}}(t_1,\cdots,t_N)
=\langle
(1+t_1f_{i_1})t_1^{-h_{i_1}}\cdots (1+t_Nf_{i_N})t_N^{-h_{i_N}}
u_{s_i\Lm_i},u_{w_0\Lm_i}\rangle.
\]
We expand $(1+t_1f_{i_1})t_1^{-h_{i_1}}\cdots (1+t_Nf_{i_N})t_N^{-h_{i_N}}
u_{s_i\Lm_i}$ as
\[
(1+t_1f_{i_1})t_1^{-h_{i_1}}\cdots (1+t_Nf_{i_N})t_N^{-h_{i_N}}
u_{s_i\Lm_i}
=\sum_{\epsilon_1,\cdots,\epsilon_N=0,1}
(t_1f_{i_1})^{\epsilon_1}t_1^{-h_{i_1}}\cdots 
(t_Nf_{i_N})^{\epsilon_N}t_N^{-h_{i_N}}u_{s_i\Lm_i}
\]
and put $F(\epsilon_1,\cdots,\epsilon_N):=(t_1f_{i_1})^{\epsilon_1}t_1^{-h_{i_1}}\cdots 
(t_Nf_{i_N})^{\epsilon_N}t_N^{-h_{i_N}}$. To prove our claim, supposing
$F(\epsilon_1,\cdots,\epsilon_N)u_{s_i\Lambda_i}=F(\epsilon_1',\cdots,\epsilon_N')u_{s_i\Lambda_i}\neq0$
for any $t_1,\cdots,t_N\in\mathbb{C}^{\times}$,
we need to show $(\epsilon_1,\cdots,\epsilon_N)=(\epsilon_1',\cdots,\epsilon_N')$.
If not, one can take $l\in[1,N]$ such that $\epsilon_l\neq \epsilon_{l'}$ and
$\epsilon_{L}=\epsilon_L'$ for $l<L\leq N$. We may assume $\epsilon_l=1$ and $\epsilon_l'=0$.
By $\epsilon_{L}=\epsilon_L'$ for $l<L\leq N$,
it can be written
\begin{eqnarray*}
& &(t_{l+1}f_{i_{l+1}})^{\epsilon_{l+1}}t_{l+1}^{-h_{i_{l+1}}}\cdots 
(t_Nf_{i_N})^{\epsilon_N}t_N^{-h_{i_N}}u_{s_i\Lambda_i}\\
&=&
(t_{l+1}f_{i_{l+1}})^{\epsilon_{l+1}'}t_{l+1}^{-h_{i_{l+1}}}\cdots 
(t_Nf_{i_N})^{\epsilon_N'}t_N^{-h_{i_N}}u_{s_i\Lambda_i}=
a u_{\mu}
\end{eqnarray*}
with an extremal weight $\mu$ and non-zero coefficient $a\in\mathbb{C}[t_{l+1}^{\pm1},\cdots,t_N^{\pm1}]$.
Since $(t_{l}f_{i_{l}})t_{l}^{-h_{i_{l}}}u_{\mu}\neq0$ by $F(\epsilon_1,\cdots,\epsilon_N)u_{s_i\Lambda_i}\neq0$,
we see that $\langle h_{i_{l}} ,\mu \rangle=1$, which yields
\begin{equation}\label{coef-pr-1}
(t_{l}f_{i_{l}})t_{l}^{-h_{i_{l}}}u_{\mu}=u_{s_{i_l}\mu}.
\end{equation}
Using $\langle h_{i_{l}} ,\mu \rangle=1$ again, it holds
\begin{equation}\label{coef-pr-2}
t_{l}^{-h_{i_{l}}}u_{\mu}=t_l^{-1}u_{\mu}.
\end{equation}
Combining (\ref{coef-pr-1}) with (\ref{coef-pr-2}), 
the coefficient of $F(\epsilon_1,\cdots,\epsilon_N)u_{s_i\Lambda_i}$ does not include the factor $t_l^{-1}$ and 
the coefficient of 
$F(\epsilon_1',\cdots,\epsilon_N')u_{s_i\Lambda_i}$ includes it, which is absurd. Thus,
we get $(\epsilon_1,\cdots,\epsilon_N)=(\epsilon_1',\cdots,\epsilon_N')$ and our claim is proved. \qed

\begin{ex}

We consider the case $\mathfrak{g}$ is of type ${\rm A}_4$ and $\textbf{i}=(4,3,2,3,1,2,4,3,2,4)$.
Following Theorem \ref{thm1}, let us compute $\Delta_{w_0\Lm_i,s_i\Lm_i}\circ \theta^-_{\textbf{i}}(t_1,\cdots,t_{10})$
for $i\in I=\{1,2,3,4\}$. First, we assume $i=1$. Since $s_4s_2s_3s_4s_2s_1s_3s_2s_3\alpha_4=\alpha_1$, the graph
$DG_0$ is the graph which has only one vertex $t_1$. The monomial in (3) of Theorem \ref{thm1} is
\[
\frac{t_5}{t_6t_9}.
\]
Since $t_1$ is the unique sink in $DG_0$, the graph $DG_1$ is as follows:
\[
DG_1:\ \ t_1\rightarrow t_1\cdot A_1^{-1}=\frac{t_2t_4}{t_7}.
\]
In $DG_1$, the monomial $\frac{t_2t_4}{t_7}$ is the unique sink and $t_4$ has a positive exponent in $\frac{t_2t_4}{t_7}$ and $t_6=t_{4^+}$ does not.
Hence, the graph $DG_2$ is as follows:
\[
DG_2:\ \ t_1\rightarrow \frac{t_2t_4}{t_7}\rightarrow \frac{t_2t_4}{t_7}\cdot A_4^{-1}=\frac{t_2t_6}{t_8}.
\]
Remark that $t_2$ has a positive exponent in $\frac{t_2t_4}{t_7}$ but $t_4=t_{2^+}$ also has a positive exponent so that we do not add the monomial
$\frac{t_2t_4}{t_7}\cdot A_2^{-1}$. We see that in $DG_2$, the monomial $\frac{t_2t_6}{t_8}$ is the unique sink
and the exponents of $t_2$, $t_6$ are positive and those of $t_{2^+}=t_4$, $t_{6^+}=t_9$ are not. Thus,
adding $\frac{t_2t_6}{t_8}\cdot A_2^{-1}=\frac{t_3t_6}{t_4t_8}$ and $\frac{t_2t_6}{t_8}\cdot A_6^{-1}=\frac{t_2}{t_9}$,
 one obtain the graph $DG_3$:
\[
\begin{xy}
(-20,0) *{DG_3:}="DG",
(0,0) *{t_1}="1",
(10,0)*{\frac{t_2t_4}{t_7}}="2",
(20,0)*{\frac{t_2t_6}{t_8}}="3",
(30,0)*{\frac{t_3t_6}{t_4t_8}}="4-1",
(20,-10)*{\frac{t_2}{t_9}}="4-2",
\ar@{->} "1";"2"^{}
\ar@{->} "2";"3"^{}
\ar@{->} "3";"4-1"^{}
\ar@{->} "3";"4-2"^{}
\end{xy}
\]
In $DG_3$, the monomials $\frac{t_3t_6}{t_4t_8}$ and $\frac{t_2}{t_9}$ are sinks. Hence,
we add the monomials $\frac{t_3t_6}{t_4t_8}\cdot A_6^{-1}=\frac{t_2}{t_9}\cdot A_2^{-1}=\frac{t_3}{t_4t_9}$:
\[
\begin{xy}
(-20,0) *{DG_4:}="DG",
(0,0) *{t_1}="1",
(10,0)*{\frac{t_2t_4}{t_7}}="2",
(20,0)*{\frac{t_2t_6}{t_8}}="3",
(30,0)*{\frac{t_3t_6}{t_4t_8}}="4-1",
(20,-10)*{\frac{t_2}{t_9}}="4-2",
(30,-10)*{\frac{t_3}{t_4t_9}}="5",
\ar@{->} "1";"2"^{}
\ar@{->} "2";"3"^{}
\ar@{->} "3";"4-1"^{}
\ar@{->} "3";"4-2"^{}
\ar@{->} "4-1";"5"^{}
\ar@{->} "4-2";"5"^{}
\end{xy}
\]
By a remark we said after Theorem \ref{thm1}, we add only one monomial $\frac{t_3}{t_4t_9}$ in this step.
Thus, the monomial $\frac{t_3}{t_4t_9}$ is the unique sink in $DG_4$. Thus, adding
the monomial $\frac{t_3}{t_4t_9}\cdot A_3^{-1}=\frac{t_5}{t_6t_9}$, one get the graph $DG_5$. 
In $DG_5$, the unique sink $\frac{t_5}{t_6t_9}$ is the monomial in (3) of Theorem \ref{thm1} so that the algorithm is finished: 
\[
\begin{xy}
(-20,0) *{DG_5=DG:}="DG",
(0,0) *{t_1}="1",
(10,0)*{\frac{t_2t_4}{t_7}}="2",
(20,0)*{\frac{t_2t_6}{t_8}}="3",
(30,0)*{\frac{t_3t_6}{t_4t_8}}="4-1",
(20,-10)*{\frac{t_2}{t_9}}="4-2",
(30,-10)*{\frac{t_3}{t_4t_9}}="5",
(40,-10)*{\frac{t_5}{t_6t_9}}="6",
\ar@{->} "1";"2"^{1}
\ar@{->} "2";"3"^{4}
\ar@{->} "3";"4-1"^{2}
\ar@{->} "3";"4-2"^{6}
\ar@{->} "4-1";"5"^{6}
\ar@{->} "4-2";"5"^{2}
\ar@{->} "5";"6"^{5}
\end{xy}
\]
Here, for two monomials $M$, $M'$, the notation $M\overset{j}{\rightarrow} M'$ implies $M'=M\cdot A_j^{-1}$.
Hence, it holds
\[
\Delta_{w_0\Lm_1,s_1\Lm_1}\circ \theta^-_{\textbf{i}}(t_1,\cdots,t_{10})
=t_1+\frac{t_2t_4}{t_7}+\frac{t_2t_6}{t_8}+\frac{t_3t_6}{t_4t_8}+\frac{t_2}{t_9}+\frac{t_3}{t_4t_9}+\frac{t_5}{t_6t_9}.
\]
Similarly, the graphs for $i=2,3,4$ are as follows:
\[
DG\ (i=2):\ t_9,
\]
\[
\begin{xy}
(-10,0) *{DG\ (i=3):}="DG",
(10,0)*{t_4}="2",
(20,0)*{\frac{t_6t_7}{t_8}}="3",
(30,0)*{\frac{t_7}{t_9}}="4-1",
(20,-10)*{\frac{t_6}{t_{10}}}="4-2",
(30,-10)*{\frac{t_8}{t_9t_{10}}}="5",
\ar@{->} "2";"3"^{4}
\ar@{->} "3";"4-1"^{6}
\ar@{->} "3";"4-2"^{7}
\ar@{->} "4-1";"5"^{7}
\ar@{->} "4-2";"5"^{6}
\end{xy},
\]
\[
DG\ (i=4):\ t_{10}.
\]
These graphs mean
\[
\Delta_{w_0\Lm_2,s_2\Lm_2}\circ \theta^-_{\textbf{i}}(t_1,\cdots,t_{10})=t_9,\quad
\Delta_{w_0\Lm_4,s_4\Lm_4}\circ \theta^-_{\textbf{i}}(t_1,\cdots,t_{10})=t_{10},
\]
\[
\Delta_{w_0\Lm_3,s_3\Lm_3}\circ \theta^-_{\textbf{i}}(t_1,\cdots,t_{10})=
t_4+\frac{t_6t_7}{t_8}+\frac{t_7}{t_9}+\frac{t_6}{t_{10}}+\frac{t_8}{t_9t_{10}}.
\]
Under the identification $X_{*}((\mathbb{C}^{\times})^{10})=\mathbb{Z}^{10}$,
we see that
\begin{eqnarray*}
& &
{\rm Trop}(\Delta_{w_0\Lm_1,s_1\Lm_1}\circ\theta^-_{\textbf{i}})(z_1,\cdots,z_{10})\\
&=&{\rm min}(z_1,z_2+z_4-z_7,z_2+z_6-z_8,z_3+z_6-z_4-z_8,z_2-z_9,z_3-z_4-z_9,z_5-z_6-z_9),
\end{eqnarray*}
\[
{\rm Trop}(\Delta_{w_0\Lm_2,s_2\Lm_2}\circ\theta^-_{\textbf{i}})(z_1,\cdots,z_{10})
=z_9,\qquad
{\rm Trop}(\Delta_{w_0\Lm_4,s_4\Lm_4}\circ\theta^-_{\textbf{i}})(z_1,\cdots,z_{10})
=z_{10},
\]
\[
{\rm Trop}(\Delta_{w_0\Lm_3,s_3\Lm_3}\circ\theta^-_{\textbf{i}})(z_1,\cdots,z_{10})
={\rm min}(z_4,z_6+z_7-z_8,z_7-z_9,z_6-z_{10},z_8-z_9-z_{10}).
\]
Hence one obtain an explicit
form of $\tilde{B}_{\theta^-_{\textbf{i}},\Phi^{h}_{\rm BK}}$:
\[
\tilde{B}_{\theta^-_{\textbf{i}},\Phi^{h}_{\rm BK}}=
\left\{(z_1,\cdots,z_{10})\in X_{*}((\mathbb{C}^{\times})^{10})\left|
\begin{array}{l} z_1\geq0,\ z_2+z_4-z_7\geq0,\ z_2+z_6-z_8\geq0,\\
z_3+z_6-z_4-z_8\geq0,\ z_2-z_9\geq0,\\
 z_3-z_4-z_9\geq0,\ z_5-z_6-z_9\geq0,\\
 z_9\geq0,\ z_4\geq0,\ z_6+z_7-z_8\geq0,\\
 z_7-z_9\geq0,\ z_6-z_{10}\geq0,\ z_8-z_9-z_{10}\geq0,\\ z_{10}\geq0.
\end{array} \right.\right\}.
\]

\end{ex}

\begin{ex}\label{D4ex}

We consider the case $\mathfrak{g}$ is of type ${\rm D}_4$ and $\textbf{i}=(1,2,3,4,2,1,4,3,2,3,4,2)$.
Following Theorem \ref{thm1}, let us compute $\Delta_{w_0\Lm_i,s_i\Lm_i}\circ \theta^-_{\textbf{i}}(t_1,\cdots,t_{12})$
for $i=1,3,4$
since $V(\Lambda_1)$, $V(\Lambda_3)$ and $V(\Lambda_4)$ are minuscule representations. First, we set $i=1$. By $s_2s_4s_3s_2s_3s_4s_1s_2s_4s_3s_2\alpha_1=\alpha_1$, it follows $DG_0$ has
the unique vertex $t_1$. Applying Theorem \ref{thm1}, we get the following decoration graph $DG$:
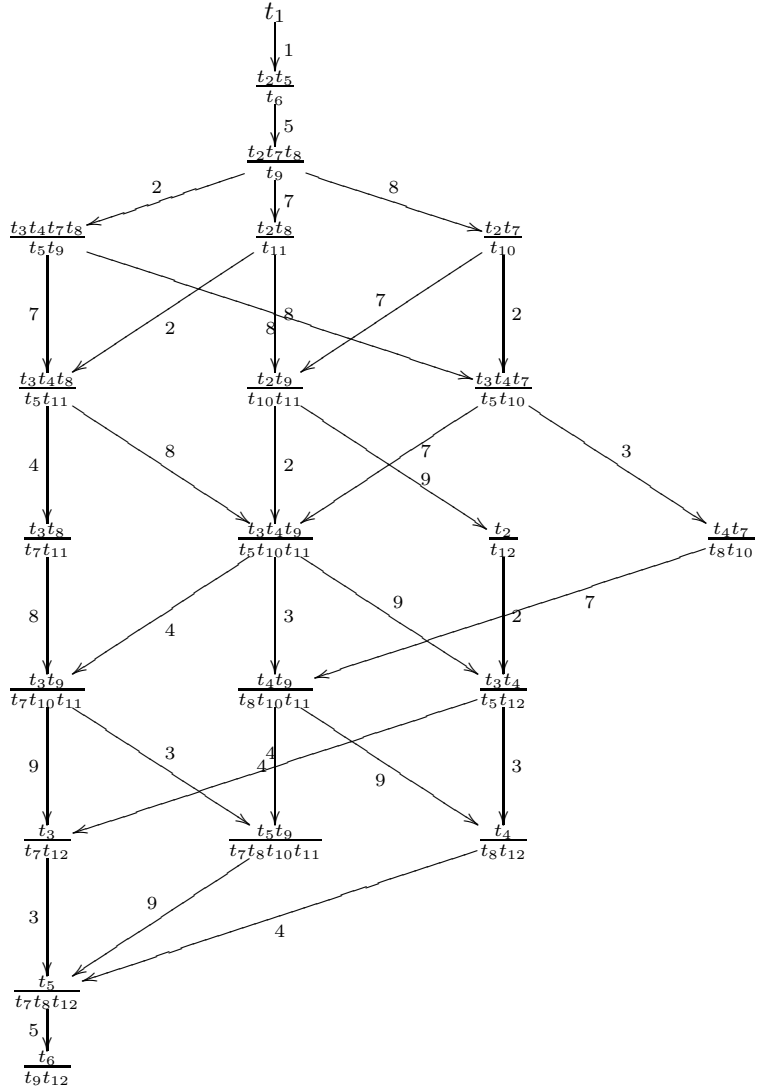
\[
\begin{xy}
(0,100) *{t_1}="1",
(0,90)*{\frac{t_2t_5}{t_6}}="2",
(0,80)*{\frac{t_2t_7t_8}{t_9}}="3",
(-30,70)*{\frac{t_3t_4t_7t_8}{t_5t_9}}="4-1",
(0,70)*{\frac{t_2t_8}{t_{11}}}="4-2",
(30,70)*{\frac{t_2t_7}{t_{10}}}="4-3",
(-30,50)*{\frac{t_3t_4t_8}{t_5t_{11}}}="5-1",
(0,50)*{\frac{t_2t_9}{t_{10}t_{11}}}="5-2",
(30,50)*{\frac{t_3t_4t_7}{t_5t_{10}}}="5-3",
(-30,30)*{\frac{t_3t_8}{t_7t_{11}}}="6-1",
(0,30)*{\frac{t_3t_4t_9}{t_5t_{10}t_{11}}}="6-2",
(30,30)*{\frac{t_2}{t_{12}}}="6-3",
(60,30)*{\frac{t_4t_7}{t_8t_{10}}}="6-4",
(-30,10)*{\frac{t_3t_9}{t_7t_{10}t_{11}}}="7-1",
(0,10)*{\frac{t_4t_9}{t_8t_{10}t_{11}}}="7-2",
(30,10)*{\frac{t_3t_4}{t_5t_{12}}}="7-3",
(-30,-10)*{\frac{t_3}{t_7t_{12}}}="8-1",
(0,-10)*{\frac{t_5t_9}{t_7t_8t_{10}t_{11}}}="8-2",
(30,-10)*{\frac{t_4}{t_8t_{12}}}="8-3",
(-30,-30)*{\frac{t_5}{t_7t_8t_{12}}}="9",
(-30,-40)*{\frac{t_6}{t_9t_{12}}}="10",
\ar@{->} "1";"2"^{1}
\ar@{->} "2";"3"^{5}
\ar@{->} "3";"4-1"_{2}
\ar@{->} "3";"4-2"^{7}
\ar@{->} "3";"4-3"^{8}
\ar@{->} "4-1";"5-1"_{7}
\ar@{->} "4-2";"5-2"^{8}
\ar@{->} "4-3";"5-3"^{2}
\ar@{->} "4-3";"5-2"_{7}
\ar@{->} "4-1";"5-3"_{8}
\ar@{->} "4-2";"5-1"^{2}
\ar@{->} "5-1";"6-1"_{4}
\ar@{->} "5-2";"6-2"^{2}
\ar@{->} "5-1";"6-2"^{8}
\ar@{->} "5-3";"6-2"_{\qquad \qquad \qquad \qquad \qquad 7}
\ar@{->} "5-2";"6-3"_{\qquad \qquad \qquad \qquad \qquad 9}
\ar@{->} "5-3";"6-4"^{3}
\ar@{->} "6-1";"7-1"_{8}
\ar@{->} "6-2";"7-1"^{4}
\ar@{->} "6-2";"7-2"^{3}
\ar@{->} "6-2";"7-3"^{9}
\ar@{->} "6-3";"7-3"^{2}
\ar@{->} "6-4";"7-2"_{\qquad \qquad \qquad \qquad \qquad 7}
\ar@{->} "7-1";"8-1"_{9}
\ar@{->} "7-1";"8-2"^{3}
\ar@{->} "7-2";"8-2"_{4}
\ar@{->} "7-2";"8-3"_{9}
\ar@{->} "7-3";"8-3"^{3}
\ar@{->} "7-3";"8-1"_{4}
\ar@{->} "8-1";"9"_{3}
\ar@{->} "8-2";"9"_{9}
\ar@{->} "8-3";"9"^{4}
\ar@{->} "9";"10"_{5}
\end{xy}
\]
Just as in the previous example, the notation $M\overset{j}{\rightarrow} M'$ implies $M'=M\cdot A_j^{-1}$.
Similarly, we obtain
\[
DG\ (i=3) :\ \ t_7\overset{7}{\rightarrow}\frac{t_9}{t_{11}}\overset{9}{\rightarrow} \frac{t_{10}}{t_{12}},\qquad
DG\ (i=4) :\ \ t_8\overset{8}{\rightarrow}\frac{t_9}{t_{10}}\overset{9}{\rightarrow} \frac{t_{11}}{t_{12}}.
\]
Although we can not apply Theorem \ref{thm1} to $\Delta_{w_0\Lm_2,s_2\Lm_2}$ since $V(\Lambda_2)$ is not minuscule,
it is easy to calculate $\Delta_{w_0\Lm_2,s_2\Lm_2}\circ \theta^-_{\textbf{i}}(t_1,\cdots,t_{12})=t_{12}$ by $i_{12}=2$
and Proposition 4.11 of \cite{BZ}. In this way, we get the explicit form of $\tilde{B}_{\theta^-_{\textbf{i}},\Phi^{h}_{\rm BK}}$
under the identification $X_{*}((\mathbb{C}^{\times})^{12})=\mathbb{Z}^{12}$ :
\begin{eqnarray*}
& &
\tilde{B}_{\theta^-_{\textbf{i}},\Phi^{h}_{\rm BK}}\\
&=&
\left\{(z_1,\cdots,z_{12})\in X_{*}((\mathbb{C}^{\times})^{12}) \left|
\begin{array}{l} z_1\geq0,\ z_2+z_5-z_6\geq0,\ z_2+z_7+z_8-z_9\geq0,\\
z_3+z_4+z_7+z_8-z_5-z_9\geq0,\ z_2+z_8-z_{11}\geq0,\\
 z_2+z_7-z_{10}\geq0,\ z_3+z_4+z_8-z_5-z_{11}\geq0,\\
 z_2+z_9-z_{10}-z_{11}\geq0,\\
 z_3+z_4+z_7-z_5-z_{10}\geq0,\\ z_3+z_8-z_7-z_{11}\geq0,
 z_3+z_4+z_9-z_5-z_{10}-z_{11}\geq0,\\
  z_2-z_{12}\geq0,\ z_4+z_7-z_8-z_{10}\geq0,\\
 z_3+z_9-z_7-z_{10}-z_{11}\geq0,\\
 z_4+z_9-z_8-z_{10}-z_{11}\geq0,\ z_3+z_4-z_5-z_{12}\geq0,\\
 z_3-z_7-z_{12}\geq0,\ z_5+z_9-z_7-z_8-z_{10}-z_{11}\geq0,\\
 z_4-z_8-z_{12}\geq0, z_5-z_7-z_8-z_{12}\geq0,\\
 z_6-z_9-z_{12}\geq0, z_{12}\geq0,\\
 z_7\geq0,\ z_9-z_{11}\geq0,\ z_{10}-z_{12}\geq0,\\
 z_8\geq0,\ z_9-z_{10}\geq0,\ z_{11}-z_{12}\geq0
\end{array} \right. \right\}.
\end{eqnarray*}

\end{ex}

\subsection{Proof of Theorem \ref{thm1}}

In what follows, $v_{-w\Lambda_i}$ denotes an extremal weight vector in $V(-w_0\Lambda_i)$ with weight $-w\Lambda_i$ for $w\in W$.
In this subsection, we give a proof of Theorem \ref{thm1}.
Let us start with the following Lemma:

\begin{lem}\label{lem1}
\begin{enumerate}
\item[(1)]
For $j\in[0,N]$ and $c_{j+1},\cdots,c_N\in\mathbb{Z}_{\geq0}$, we suppose that
$e^{c_{j+1}}_{i_{j+1}}\cdots e^{c_N}_{i_N}v_{-s_i\Lambda_i}\neq 0$.
For $t\in I$,
the value of
\[
\langle h_t,\ {\rm wt}(e^{c_{j+1}}_{i_{j+1}}\cdots e^{c_N}_{i_N}v_{-s_i\Lambda_i}) \rangle
\]
is $1$ or $0$ or $-1$.
\item[(2)] Let
$\pi=(\gamma_0,\gamma_1,\cdots,\gamma_N)$ be an \textbf{i}-trail
from $-w_0\Lambda_i$ to $-s_i\Lambda_i$. For $j\in[1,N]$ and $t\in I$, the values of
\[
\gamma_{j-1}(h_t),\ \gamma_{j}(h_t)
\]
\[
d_j(\pi)=\frac{\gamma_{j-1}+\gamma_j}{2}(h_{i_j})
\]
are $1$ or $0$ or $-1$.
\end{enumerate}
\end{lem}

\nd
{\it Proof.}

\nd
(1) Our claim
is an easy consequence of (\ref{mini-pro}).

\vspace{2mm}

\nd
(2) Let $\{c_l\}$ be the integers for $\pi$ as in Definition \ref{pretrail}.
Since
\[
\gamma_j
=\gamma_N+\sum^{N-1}_{t=j}(\gamma_t-\gamma_{t+1}) 
= -s_i\Lambda_i + \sum^{N-1}_{t=j} c_{t+1}\alpha_{i_{t+1}} 
= {\rm wt}(e^{c_{j+1}}_{i_{j+1}}\cdots e^{c_N}_{i_N}v_{-s_i\Lambda_i}),
\]
\[
\gamma_{j-1}
= {\rm wt}(e^{c_{j}}_{i_{j}} e^{c_{j+1}}_{i_{j+1}}\cdots e^{c_N}_{i_N}v_{-s_i\Lambda_i}),
\]
our statement (2) follows by (1). \qed

\begin{lem}\label{lem2}
If  a monomial $M$ appearing in
$\Delta_{w_0\Lm_i,s_i\Lm_i}\circ \theta^-_{\textbf{i}}(t_1,\cdots,t_N)$
belongs to $\{ \prod^N_{l=1} t_l^{d_l} | d_l\in\mathbb{Z}_{\geq0} \}$ then
\[
M=t_k,
\]
where $k\in[1,N]$ is the same one as in Theorem \ref{thm1}. 
\end{lem}

\nd
{\it Proof.}

Let $\pi=(\gamma_0,\gamma_1,\cdots,\gamma_N)$ be the \textbf{i}-trail
from $-w_0\Lambda_i$ to $-s_i\Lambda_i$ corresponding to $M$, that is, putting
\[
d_l=\frac{\gamma_{l-1}+\gamma_l}{2}(h_{i_l}),
\]
it holds $M=\prod^N_{l=1} t_l^{d_l}$ (Theorem \ref{trail-thm}).
As seen in the proof of Lemma \ref{lem1} (2), for each $l\in[1,N]$, it holds
\begin{equation}\label{pr-3}
\gamma_l
= {\rm wt}(e^{c_{l+1}}_{i_{l+1}}\cdots e^{c_N}_{i_N}v_{-s_i\Lambda_i}),\qquad
\gamma_{l-1}
= {\rm wt}(e^{c_{l}}_{i_{l}} e^{c_{l+1}}_{i_{l+1}}\cdots e^{c_N}_{i_N}v_{-s_i\Lambda_i}),
\end{equation}
where the integers $\{c_l\}$ are same as in Definition \ref{pretrail}.
Hence we see that $\gamma_l(h_{i_l})\leq \gamma_l(h_{i_l})+2c_l=\gamma_{l-1}(h_{i_l})$.
The nonnegativity of $d_l$ and Lemma \ref{lem1} (2) mean that $\gamma_l(h_{i_l})=\gamma_{l-1}(h_{i_l})=1$ or 
$\gamma_l(h_{i_l})=\gamma_{l-1}(h_{i_l})=0$ or $\gamma_l(h_{i_l})=-1$, $\gamma_{l-1}(h_{i_l})=1$.
Hence, if $\gamma_l(h_{i_l})=0$ (resp. $\gamma_l(h_{i_l})=1$,\ $\gamma_l(h_{i_l})=-1$)
then $\gamma_{l-1}(h_{i_l})=0$ (resp. $\gamma_{l-1}(h_{i_l})=1$, $\gamma_{l-1}(h_{i_l})=1$). Thus,
the value of $\gamma_{l-1}(h_{i_l})$ and
\[
c_l=\frac{\gamma_{l-1}-\gamma_{l}}{2}(h_{i_l})
\]
are uniquely determined from $\gamma_{l}(h_{i_l})$.
In particular,
since $\gamma_{N}(h_{i_N})=-s_i\Lambda_i(h_{i_N})$ is a fixed value, the value of $c_N$ is
uniquely determined. By the relation in (\ref{pr-3}), the weight $\gamma_{N-1}$
is also
uniquely determined. From the value of $\gamma_{N-1}(h_{i_{N-1}})$, the integer
$c_{N-1}$ is
uniquely determined
and so is weight $\gamma_{N-2}$ by (\ref{pr-3}). Repeating this argument, we see that $\gamma_0,\gamma_1,\gamma_2,\cdots,\gamma_N$ are uniquely determined and so is the $\textbf{i}$-trail $\pi$.
Since Proposition \ref{high-prop} says $t_k\in \{ \prod^N_{l=1} t_l^{d_l} | d_l\in\mathbb{Z}_{\geq0}\}$ appears
in $\Delta_{w_0\Lm_i,s_i\Lm_i}\circ \theta^-_{\textbf{i}}(t_1,\cdots,t_N)$ as a term, we get $M=t_k$. \qed

\begin{lem}\label{jpm}
We suppose $j\in[1,N]$ and $j^+\leq N$.
For a term $t_1^{d_1}\cdots t_N^{d_N}$ of $\Delta_{w_0\Lm_i,s_i\Lm_i}\circ \theta^-_{\textbf{i}}(t_1,\cdots,t_N)$, if $d_j=1$ then
$d_{j^+}\neq -1$.
\end{lem}
\nd
{\it Proof.}

Let $\pi=(\gamma_0,\gamma_1,\cdots,\gamma_N)$ be the $\textbf{i}$-trail corresponding 
to the term $t_1^{d_1}\cdots t_N^{d_N}$ with integers $\{c_l\}$ in Definition \ref{pretrail}. 
The assumption $d_j=\frac{\gamma_{j-1}+\gamma_j}{2}(h_{i_j})=1$ and Lemma \ref{lem1} (2)
yields $\gamma_j(h_{i_j})=\gamma_{j-1}(h_{i_j})=1$.
The equation $\gamma_j(h_{i_j})=1$ means
\[
\langle h_{i_j},
{\rm wt}(e^{c_{j+1}}_{i_{j+1}}\cdots e^{c_N}_{i_N}v_{-s_i\Lambda_i})
\rangle=1.
\]
It is easy to see
\begin{eqnarray*}
\gamma_{j^+-1}(h_{i_j})&=&\langle h_{i_j},
{\rm wt}(e^{c_{j^+}}_{i_{j^+}}\cdots e^{c_N}_{i_N}v_{-s_i\Lambda_i})
\rangle\\
&=&\langle h_{i_j},
{\rm wt}(e^{c_{j+1}}_{i_{j+1}}\cdots e^{c_{j^+}}_{i_{j^+}}\cdots e^{c_N}_{i_N}v_{-s_i\Lambda_i})
\rangle
-\langle h_{i_j},\sum^{j^+-1}_{t=j+1} c_t\alpha_{i_t} \rangle\geq1.
\end{eqnarray*}
Considering Lemma \ref{lem1} (2), we obtain $\gamma_{j^+-1}(h_{i_j})=1$.
Since $d_{j^+}=\frac{\gamma_{j^+-1}+\gamma_{j^+}}{2}(h_{i_j})$, 
using Lemma \ref{lem1} (2) again, it holds $d_{j^+}\geq0$. \qed

\begin{prop}\label{term-prop}
We suppose $j\in[1,N]$ and $j^+\leq N$.
For a term $t_1^{d_1}\cdots t_N^{d_N}$ of $\Delta_{w_0\Lm_i,s_i\Lm_i}\circ \theta^-_{\textbf{i}}(t_1,\cdots,t_N)$, it holds
$d_j=1$ and $d_{j^+}=0$ if and only if the corresponding $\textbf{i}$-trail satisfies
\begin{equation}\label{cond-1}
c_j=0,\ c_{j^+}=1,\ \text{for }s\in[j+1,j^+-1]\ \text{such that }a_{i_s,i_j}<0,\ \text{it holds }c_s=0.
\end{equation}
Here, $\{c_l\}$ are integers in Definition \ref{pretrail}. 
\end{prop}
\nd
{\it Proof.}

Let $\pi=(\gamma_0,\gamma_1,\cdots,\gamma_N)$ be the $\textbf{i}$-trail corresponding 
to the term $t_1^{d_1}\cdots t_N^{d_N}$ with integers $\{c_l\}$ in Definition \ref{pretrail}.
First, we suppose that $d_j=1$ and $d_{j^+}=0$.
By $d_j=\frac{\gamma_{j-1}+\gamma_j}{2}(h_{i_j})=1$ and Lemma \ref{lem1} (2), it holds $\gamma_j(h_{i_j})=\gamma_{j-1}(h_{i_j})=1$.
It follows from $\gamma_{j-1}-\gamma_j=c_j\alpha_{i_j}$
that $c_j=0$. The equality $\gamma_j(h_{i_j})=1$ implies
\begin{equation}\label{pr-0}
\langle h_{i_j},
{\rm wt}(e^{c_{j+1}}_{i_{j+1}}\cdots e^{c_N}_{i_N}v_{-s_i\Lambda_i})
\rangle=1.
\end{equation}
Let us take $s\in[j+1,j^+-1]$ such that $a_{i_s,i_j}<0$. If $c_s>0$ then
(\ref{pr-0}) means
\begin{eqnarray*}
& &\langle h_{i_j},
{\rm wt}(e^{c_{s+1}}_{i_{s+1}}\cdots e^{c_N}_{i_N}v_{-s_i\Lambda_i})
\rangle\\
&=&\langle h_{i_j},
{\rm wt}(e^{c_{j+1}}_{i_{j+1}}\cdots e^{c_{s+1}}_{i_{s+1}}\cdots e^{c_N}_{i_N}v_{-s_i\Lambda_i})
\rangle
-\langle h_{i_j},\sum^{s}_{t=j+1} c_t\alpha_{i_t} \rangle
>1,
\end{eqnarray*}
which contradicts Lemma \ref{lem1} (1).
Therefore, it holds $c_s=0$. In the above formula, taking $s=j^+-1$ and
combining with (\ref{pr-0}), it follows
\[
\gamma_{j^+ -1}(h_{i_j})
=
\langle h_{i_j},
{\rm wt}(e^{c_{j^+}}_{i_{j^+}}\cdots e^{c_N}_{i_N}v_{-s_i\Lambda_i})
\rangle=1.
\]
By $d_{j^+}=0$, it holds $\gamma_{j^+}(h_{i_j})+\gamma_{j^+-1}(h_{i_j})=0$ so that
$\gamma_{j^+}(h_{i_j})=-1$. It follows by $\gamma_{j^+-1}-\gamma_{j^+}=c_{j^+}\alpha_{i_j}$ that
$c_{j^+}=1$. Hence, the condition (\ref{cond-1}) holds.

Next, we suppose the corresponding $\textbf{i}$-trail satisfies (\ref{cond-1}).
The assumption $c_{j^+}=\frac{\gamma_{j^+-1}-\gamma_{j^+}}{2}(h_{i_j})=1$ and Lemma \ref{lem1} (2) imply that 
$\gamma_{j^+-1}(h_{i_j})=1$ and $\gamma_{j^+}(h_{i_j})=-1$.
Therefore, we obtain $d_{j^+}=\frac{\gamma_{j^+-1}+\gamma_{j^+}}{2}(h_{i_j})=0$.
We see that
\[
\gamma_{j^+-1}(h_{i_j})
=\langle h_{i_j},
{\rm wt}(e^{c_{j^+}}_{i_{j^+}}
e^{c_{j^++1}}_{i_{j^++1}}\cdots e^{c_N}_{i_N}v_{-s_i\Lambda_i})
\rangle =1.
\]
Since we supposed for $s\in[j+1,j^+-1]$ such that
$a_{i_s,i_j}<0$, it follows $c_s=0$, we can verify that
\[
\gamma_{j}(h_{i_j})=
\langle h_{i_j},
{\rm wt}(e^{c_{j+1}}_{i_{j+1}}\cdots e^{c_{j^+}}_{i_{j^+}}\cdots e^{c_N}_{i_N}v_{-s_i\Lambda_i})
\rangle =1.
\]
The assumption $c_j=\frac{\gamma_{j-1}-\gamma_{j}}{2}(h_{i_j})=0$ means that
$\gamma_{j-1}(h_{i_j})=1$, which yields $d_j=\frac{\gamma_{j-1}+\gamma_{j}}{2}(h_{i_j})=1$. \qed

\begin{prop}\label{term-prop2}
Let $j\in[1,N]$.
For a term $t_1^{d_1}\cdots t_N^{d_N}$ of $\Delta_{w_0\Lm_i,s_i\Lm_i}\circ \theta^-_{\textbf{i}}(t_1,\cdots,t_N)$,
if $d_j=1$ and $j^+=N+1$ then $i_j=i$ and
\[
t_1^{d_1}\cdots t_N^{d_N}
=t_jt_{j+1}^{a_{i_{j+1},i_j}}\cdots t_N^{^{a_{i_{N},i_j}}}.
\]
\end{prop}
\nd
{\it Proof.}

Let $\pi=(\gamma_0,\gamma_1,\cdots,\gamma_N)$ be the $\textbf{i}$-trail corresponding 
to the term $t_1^{d_1}\cdots t_N^{d_N}$ with integers $\{c_l\}$.
By $d_j=\frac{\gamma_{j-1}+\gamma_j}{2}(h_{i_j})=1$,
it follows by Lemma \ref{lem1} (2) that $\gamma_{j-1}(h_{i_j})=\gamma_{j}(h_{i_j})=1$ so that $c_j=\frac{\gamma_{j-1}-\gamma_j}{2}(h_{i_j})=0$.

We suppose that $i_j\neq i$. Note that by $s_i\Lambda_i=\Lambda_i-\alpha_i$, we obtain
\[
\langle h_{i_j}, {\rm wt}(v_{-s_i\Lambda_i})\rangle\leq0.
\]
Thus, the assumption $j^+=N+1$ means
\[
1=\gamma_{j}(h_{i_j})=\langle h_{i_j},
{\rm wt}(e^{c_{j+1}}_{i_{j+1}}\cdots e^{c_N}_{i_N}v_{-s_i\Lambda_i})
\rangle
=\langle h_{i_j}, {\rm wt}(v_{-s_i\Lambda_i})\rangle
+ \langle h_{i_j},\sum^{N}_{t=j+1} c_t\alpha_{i_t} \rangle
 \leq0,
\]
which is absurd.
Hence, we suppose $i_j= i$. It follows 
\begin{equation}\label{pr-1}
\langle h_{i_j}, {\rm wt}(v_{-s_i\Lambda_i})\rangle=1.
\end{equation}
By (\ref{pr-1}) and
\[
1=\gamma_{j}(h_{i_j})=\langle h_{i_j},
{\rm wt}(e^{c_{j+1}}_{i_{j+1}}\cdots e^{c_N}_{i_N}v_{-s_i\Lambda_i})\rangle
=\langle h_{i_j}, {\rm wt}(v_{-s_i\Lambda_i})\rangle
+
\langle h_{i_j}, \sum^N_{t=j+1} c_t\alpha_{i_t} \rangle,
\]
for any $l\in[j+1,N]$ such that $a_{i_j,i_l}<0$, it holds $c_l=0$.
 For $l\in[j+1,N]$ such that $a_{i_j,i_l}=0$, it holds $e_{i_l}v_{-s_i\Lambda_i}=0$ by $\langle h_{i_l},-s_i\Lambda_i \rangle=0$. Thus, $c_l=0$.
In this way, we see that $c_{j+1}=c_{j+2}=\cdots =c_N=0$.
The exponents $c_1,\cdots,c_{j-1}$ are uniquely determined by
\[
c_l={\rm max}\{t\in\mathbb{Z}_{\geq0} | e^t_{i_{l}} e^{c_{l+1}}_{i_{l+1}}\cdots e^{c_{j-1}}_{i_{j-1}}v_{-s_i\Lambda_i} \}
\]
for $l=1,2,\cdots,j-1$. Hence, the $\textbf{i}$-trail is uniquely determined.
Since the function $\Delta_{w_0\Lm_i,s_i\Lm_i}\circ \theta^-_{\textbf{i}}(t_1,\cdots,t_N)$
has a term (\ref{lowest-term}) and this term satisfies the assumption of our claim,
it holds
\[
t_1^{d_1}\cdots t_N^{d_N}
=t_jt_{j+1}^{a_{i_{j+1},i_j}}\cdots t_N^{^{a_{i_{N},i_j}}}.
\]
\qed

\begin{prop}\label{term-prop2-a}
For a term $t_1^{d_1}\cdots t_N^{d_N}$ of $\Delta_{w_0\Lm_i,s_i\Lm_i}\circ \theta^-_{\textbf{i}}(t_1,\cdots,t_N)$,
it holds $d_j>0$ for some $j\in[1,N]$.
\end{prop}
\nd
{\it Proof.}

We assume $d_j\leq0$ for all $j\in[1,N]$ and deduce a contradiction from this assumption.
Let $\pi=(\gamma_0,\gamma_1,\cdots,\gamma_N)$ be the $\textbf{i}$-trail corresponding 
to the term $t_1^{d_1}\cdots t_N^{d_N}$ with integers $\{c_l\}$.
One define $i^*\in I$ as the index such that $\gamma_0=-w_0\Lambda_i=\Lambda_{i^*}$.
Let us prove
\begin{equation}\label{2-a-1}
\gamma_l=s_{i_l}s_{i_{l-1}}\cdots s_{i_1}\Lambda_{i^*}\ (l\in[0,N])
\end{equation}
by the induction on $l$. We suppose $\gamma_{l-1}=s_{i_{l-1}}s_{i_{l-2}}\cdots s_{i_1}\Lambda_{i^*}$.
Since $s_{i_1}\cdots s_{i_{l-2}}s_{i_{l-1}}h_{i_l}$ is a positive coroot, it holds
\[
\gamma_{l-1}(h_{i_l})=
\langle 
h_{i_l}, s_{i_{l-1}}s_{i_{l-2}}\cdots s_{i_1}\Lambda_{i^*}
\rangle
=\langle 
s_{i_1}\cdots s_{i_{l-2}}s_{i_{l-1}}h_{i_l},\Lambda_{i^*}
\rangle\in \{0,1\}.
\]
By $d_l=\frac{\gamma_{l-1}(h_{i_l})+\gamma_{l}(h_{i_l})}{2}\leq0$, we see that
if $\gamma_{l-1}(h_{i_l})=0$ (resp. $\gamma_{l-1}(h_{i_l})=1$) then $\gamma_{l}(h_{i_l})=0$ (resp. $\gamma_{l}(h_{i_l})=-1$)
and $c_l=\frac{\gamma_{l-1}(h_{i_l})-\gamma_{l}(h_{i_l})}{2}=0$ (resp. $c_l=1$). Thus,
it follows $\gamma_l=\gamma_{l-1}-c_l\alpha_{i_l}=\gamma_{l-1}-\gamma_{l-1}(h_{i_l})\alpha_{i_l}=
s_{i_l}\gamma_{l-1}$, which induces (\ref{2-a-1}). In particular, one gets $\gamma_N=w_0\Lambda_{i^*}=-\Lambda_i$,
which contradicts $\gamma_N=-s_i\Lambda_i$. \qed

\begin{prop}\label{term-prop3}
We suppose $j\in[1,N]$ and $j^-\geq 1$.
For a term $t_1^{d_1}\cdots t_N^{d_N}$ of $\Delta_{w_0\Lm_i,s_i\Lm_i}\circ \theta^-_{\textbf{i}}(t_1,\cdots,t_N)$, it holds
$d_j=-1$ and $d_{j^-}=0$ if and only if the corresponding $\textbf{i}$-trail satisfies
\begin{equation}\label{cond-2}
c_j=0,\ c_{j^-}=1,\ \text{for }s\in[j^- +1,j-1]\ \text{such that }a_{i_s,i_j}<0,\ \text{it holds }c_s=0.
\end{equation}
\end{prop}
\nd
{\it Proof.}

Let $\pi=(\gamma_0,\gamma_1,\cdots,\gamma_N)$ be the $\textbf{i}$-trail with integers $\{c_l\}$ corresponding 
to the term $t_1^{d_1}\cdots t_N^{d_N}$. First, we suppose that
$d_j=-1$ and $d_{j^-}=0$. By $d_j=-1$, it holds
\[
d_j=\frac{\gamma_{j-1}+\gamma_j}{2}(h_{i_j})=-1
\]
and $\gamma_{j-1}(h_{i_j})=\gamma_{j}(h_{i_j})=-1$ by Lemma \ref{lem1} (2).
Therefore, we have $c_j=\frac{\gamma_{j-1}(h_{i_j})-\gamma_{j}(h_{i_j})}{2}=0$.
We take $s\in[j^- +1,j-1]$ such that $a_{i_s,i_j}<0$. 
Since 
\[
\gamma_{j-1}(h_{i_j})=\langle h_{i_j}, {\rm wt}(e^{c_{j}}_{i_{j}}\cdots e^{c_N}_{i_N}v_{-s_i\Lambda_i}) \rangle =-1,
\]
if $c_s>0$ then
\[
\langle h_{i_j}, {\rm wt}(e^{c_{s}}_{i_{s}}\cdots e^{c_{j}}_{i_{j}}\cdots e^{c_N}_{i_N}v_{-s_i\Lambda_i}) \rangle =
\langle h_{i_j}, {\rm wt}(e^{c_{j}}_{i_{j}}\cdots e^{c_N}_{i_N}v_{-s_i\Lambda_i}) \rangle
+\langle h_{i_j}, \sum^{j-1}_{t=s} c_t\alpha_{i_t} \rangle
<-1,
\]
which contradicts Lemma \ref{lem1} (1). Thus, $c_s=0$ and 
\begin{equation}\label{pr-2}
\gamma_{j^-}(h_{i_j})
=\langle h_{i_j}, {\rm wt}(e^{c_{j^-+1}}_{i_{j^-+1}}\cdots e^{c_{j-1}}_{i_{j-1}} e^{c_{j}}_{i_{j}}\cdots e^{c_N}_{i_N}v_{-s_i\Lambda_i}) \rangle
=-1.
\end{equation}
It follows from the assumption $d_{j^-}=\frac{\gamma_{j^- -1}+\gamma_{j^-}}{2}(h_{i_j})=0$ that
$\gamma_{j^- -1}(h_{i_j})=1$.
In conjunction with (\ref{pr-2}), we get $c_{j^-}=\frac{\gamma_{j^--1}-\gamma_{j^-}}{2}(h_{i_j})=1$ and the $\textbf{i}$-trail satisfies (\ref{cond-2}).

Next, we assume (\ref{cond-2}). Taking the assumption $c_{j^-}=\frac{\gamma_{j^--1}(h_{i_j})-\gamma_{j^-}(h_{i_j})}{2}=1$ into account, we see that $\gamma_{j^--1}(h_{i_j})=1$ and $\gamma_{j^-}(h_{i_j})=-1$ by Lemma \ref{lem1} (2)
so that
\[
d_{j^-}=
\frac{\gamma_{j^--1}+\gamma_{j^-}}{2}(h_{i_j})=0.
\]
By $\gamma_{j^-}(h_{i_j})=-1$, we obtain
\[
\langle h_{i_j}, {\rm wt}(e^{c_{j^-+1}}_{i_{j^-+1}}\cdots e^{c_N}_{i_N}v_{-s_i\Lambda_i})
\rangle =\gamma_{j^-}(h_{i_j})=-1.
\]
Since $c_s=0$
for
$s\in[j^- +1,j-1]$ such that $a_{i_s,i_j}<0$, it holds
\[
\gamma_{j-1}(h_{i_j})
=\langle h_{i_j}, {\rm wt}(e^{c_{j}}_{i_{j}}\cdots e^{c_N}_{i_N}v_{-s_i\Lambda_i})
\rangle =-1.
\]
Using the assumption $c_j=\frac{\gamma_{j-1}-\gamma_{j}}{2}(h_{i_j})=0$, it follows
$\gamma_{j}(h_{i_j})=-1$.
Consequently,
\[
d_j=\frac{\gamma_{j-1}+\gamma_j}{2}(h_{i_j})=-1.
\]
\qed

\begin{prop}\label{term-prop4}
Let $j\in[1,N]$.
For a term $t_1^{d_1}\cdots t_N^{d_N}$ of $\Delta_{w_0\Lm_i,s_i\Lm_i}\circ \theta^-_{\textbf{i}}(t_1,\cdots,t_N)$,
if $d_j=-1$ then $j^-\geq1$.
\end{prop}

\nd
{\it Proof.}

Let $\pi=(\gamma_0,\gamma_1,\cdots,\gamma_N)$ be the $\textbf{i}$-trail corresponding 
to the term $t_1^{d_1}\cdots t_N^{d_N}$. There exist positive integers $c_l$ ($l=1,2,\cdots,N$) such
that $\gamma_{l-1}-\gamma_l=c_l\alpha_{i_l}$. Considering Lemma \ref{lem1} (2),
the assumption
$d_j=-1$ means
\[
d_j=\frac{\gamma_{j-1}+\gamma_j}{2}(h_{i_j})=-1
\]
and $\gamma_{j-1}(h_{i_j})=\gamma_{j}(h_{i_j})=-1$. By
\[
-1=\gamma_{j-1}(h_{i_j})=\langle h_{i_j}, {\rm wt}(e^{c_{j}}_{i_{j}}\cdots e^{c_N}_{i_N}v_{-s_i\Lambda_i})
\rangle,
\]
if $j^-=0$ then
\[
\langle h_{i_j}, {\rm wt}(v_{-w_0\Lambda_i}) \rangle=
\langle h_{i_j}, {\rm wt}(e^{c_{1}}_{i_{1}}\cdots e^{c_{j}}_{i_{j}}\cdots e^{c_N}_{i_N}v_{-s_i\Lambda_i}) \rangle
\leq-1.
\]
Since $-w_0\Lambda_i$ is a dominant integral weight, it is absurd.
Consequently, we get $j^-\geq1$. \qed

\begin{prop}\label{trail-prop}
Let $\pi=(\gamma_0,\gamma_1,\cdots,\gamma_N)$, $\pi'=(\gamma_0',\gamma_1',\cdots,\gamma_N')$ be pre-\textbf{i}-trails
from $-w_0\Lambda_i$ to $-s_i\Lambda_i$. We take nonnegative integers $c_l$, $c_l'$ $(l=1,2,\cdots,N)$ as
\[
\gamma_{l-1}-\gamma_l=c_l\alpha_{i_l},\quad \gamma_{l-1}'-\gamma_l'=c_l'\alpha_{i_l}.
\]
We suppose that there exists $j\in[1,N]$ such that $j^+\leq N$ and
\[
c_l'=c_{l}\quad \text{for}\ l\in[1,N]\setminus \{j,j^+\},
\]
\[
c_j'=c_{j}+1,\quad c_{j^+}'=c_{j^+}-1.
\]
Then
\[
t_1^{d_1(\pi')}\cdots t_N^{d_N(\pi')}=t_1^{d_1(\pi)}\cdots t_N^{d_N(\pi)}A^{-1}_{j},
\]
where $d_m$ is defined as (\ref{dk}) and $A_j$ is defined in (\ref{ajdef}).
\end{prop}

\nd
{\it Proof.}

The assumption implies that for $l\in[1,N]\setminus [j,j^+-1]$,
\[
\gamma_l'=\gamma_l
\]
by $\gamma_l=\gamma_{l+1}+c_{l+1}\alpha_{i_{l+1}}$ and $\gamma_N=\gamma'_N$.
For $l\in[j,j^+-1]$, it holds
\[
\gamma_l'=\gamma_l-\alpha_{i_j}.
\]
Hence, for $l\in[1,N]\setminus [j,j^+]$,
\[
d_l(\pi')=
\frac{\gamma_{l-1}'+\gamma_l'}{2}(h_{i_l})
=
\frac{\gamma_{l-1}+\gamma_l}{2}(h_{i_l})=d_l(\pi).
\]
We also see that
\[
d_j(\pi')=
\frac{\gamma_{j-1}'+\gamma_j'}{2}(h_{i_j})
=
\frac{\gamma_{j-1}+\gamma_j-\alpha_{i_j}}{2}(h_{i_j})
=d_j(\pi)-1,
\]
\[
d_{j^+}(\pi')=
\frac{\gamma_{j^+-1}'+\gamma_{j^+}'}{2}(h_{i_j})
=
\frac{\gamma_{j^+-1}-\alpha_{i_j}+\gamma_{j^+}}{2}(h_{i_j})
=d_{j^+}(\pi)-1
\]
and for $l\in[j+1,j^+-1]$,
\[
d_l(\pi')=
\frac{\gamma_{l-1}'+\gamma_l'}{2}(h_{i_l})
=
\frac{\gamma_{l-1}+\gamma_l-2\alpha_{i_j}}{2}(h_{i_l})
=d_l(\pi)-a_{i_l,i_j}.
\]
Therefore, it holds
\[
t_1^{d_1(\pi')}\cdots t_N^{d_N(\pi')}=t_1^{d_1(\pi)}\cdots t_N^{d_N(\pi)}A^{-1}_{j}.
\]
\qed

Let $M=t_1^{d_1}\cdots t_N^{d_N}$ be 
a monomial appearing in $\Delta_{w_0\Lm_i,s_i\Lm_i}\circ \theta^-_{\textbf{i}}(t_1,\cdots,t_N)$
and
$\pi=(\gamma_0,\gamma_1,\cdots,\gamma_N)$
be the $\textbf{i}$-trail from $-w_0\Lambda_i$ to $-s_i\Lambda_i$
corresponding to $M$ with nonnegative integers $c_l=\frac{\gamma_{l-1}-\gamma_l}{2}(h_{i_l})$ $(l\in[1,N])$.
We set
\begin{equation}\label{L-def}
L(M):=\sum^{N}_{l=1} l\cdot c_l \in\mathbb{Z}_{\geq0}.
\end{equation}

\begin{lem}\label{lem3}
Let $M$ be a monomial appearing in $\Delta_{w_0\Lm_i,s_i\Lm_i}\circ \theta^-_{\textbf{i}}(t_1,\cdots,t_N)$.
If $M\neq t_k$ then there exist $m\in[1,N]$ and a monomial $M'$ in $\Delta_{w_0\Lm_i,s_i\Lm_i}\circ \theta^-_{\textbf{i}}(t_1,\cdots,t_N)$
such that $m^+\leq N$, $M=M'\cdot A_{m}^{-1}$, $L(M)<L(M')$ and in $M'$, the exponents of $t_m$ and $t_{m^+}$ are equal to $1$ and $0$, respectively. 
\end{lem}
\nd
{\it Proof.}

Let $\pi=(\gamma_0,\gamma_1,\cdots,\gamma_N)$
be the $\textbf{i}$-trail corresponding to $M=t_1^{d_1}\cdots t_N^{d_N}$
with integers $c_l=\frac{\gamma_{l-1}-\gamma_l}{2}(h_{i_l})$ $(l\in[1,N])$.
Taking Lemma \ref{lem2} into account, we can take an integer $j\in[1,N]$ as
\begin{equation}\label{defj}
j:={\rm min}\{ l\in[1,N] | d_l=-1 \}.
\end{equation}
It holds $j^-\geq1$ by Proposition \ref{term-prop4}. Since $d_j=-1=\frac{\gamma_{j-1}+\gamma_j}{2}(h_{i_j})$,
we have $\gamma_{j-1}(h_{i_j})=\gamma_{j}(h_{i_j})=-1$. Considering
\[
-1=\gamma_{j-1}(h_{i_j})
=\langle h_{i_j}, {\rm wt}(e^{c_{j}}_{i_{j}}\cdots e^{c_N}_{i_N}v_{-s_i\Lambda_i}) \rangle,
\]
it holds
\[
\gamma_{j^-}(h_{i_j})
=\langle h_{i_j}, {\rm wt}(e^{c_{j^- +1}}_{i_{j^- +1}}\cdots e^{c_{j}}_{i_{j}}\cdots e^{c_N}_{i_N}v_{-s_i\Lambda_i}) \rangle\leq -1.
\]
Lemma \ref{lem1} (2) says $\gamma_{j^-}(h_{i_j})=-1$. 
By the minimality of $j$ (\ref{defj}), it holds
$d_{j^-}=\frac{\gamma_{j^--1}+\gamma_{j^-}}{2}(h_{i_{j}})\geq0$.
Hence, Lemma \ref{lem1} (2) and $\gamma_{j^-}(h_{i_j})=-1$ imply
$\gamma_{j^--1}(h_{i_j})=1$ and $d_{j^-}=0$.
Using Proposition \ref{term-prop3}, the $\textbf{i}$-trail $\pi$ satisfies $c_j=0$, $c_{j^-}=1$ and for $s\in[j^- +1,j-1]$
such that $a_{i_s,i_j}<0$, it follows $c_s=0$.
Let $\pi'=(\gamma_0',\gamma_1',\cdots,\gamma_N')$ be the pre-\textbf{i}-trail
with integers $c_l'=\frac{\gamma_{l-1}'-\gamma_l'}{2}(h_{i_l})$ $(l\in[1,N])$ such that
\[
c_j'=1,\ c_{j^-}'=0,\quad c_l'=c_l \ \ (l\in[1,N]\setminus \{j,j^-\}).
\]
We can verify $\pi'$ is an \textbf{i}-trail since $e_te_{t'}=e_{t'}e_t$ if $a_{t,t'}=0$. Hence,
putting $M':=t_1^{d_1(\pi')}\cdots t_N^{d_N(\pi')}$, the Laurent monomial $M'$ is a term
appearing in $\Delta_{w_0\Lm_i,s_i\Lm_i}\circ \theta^-_{\textbf{i}}(t_1,\cdots,t_N)$
and
using Proposition \ref{trail-prop}, it holds $M=M'\cdot A_{j^-}^{-1}$.
It is easy to see $d_{j^-}(\pi')=1$ and $d_{j}(\pi')=0$ using Proposition \ref{term-prop}.
Putting $m=j^-$, we get our claim.\qed

\vspace{2mm}

\nd
{\it Proof of Theorem \ref{thm1}}

Let $\mathcal{M}$ be the set of monomials 
appearing in
$\Delta_{w_0\Lm_i,s_i\Lm_i}\circ \theta^-_{\textbf{i}}(t_1,\cdots,t_N)$. For a directed graph $H$,
let $V(H)$ be the set of vertices in $H$.
First, we prove $V(DG)\subset \mathcal{M}$. By Proposition \ref{high-prop}, we see that $V(DG_0)=\{t_k\}\subset \mathcal{M}$.
We prove $V(DG_l)\subset \mathcal{M}$ by induction on $l$.
We assume that $V(DG_l)\subset \mathcal{M}$ for some $l\in\mathbb{Z}_{\geq0}$. 
Let $M\in V(DG_l)$ be a sink which has a factor $t_j^{+1}$ such that $j^+\leq N$ and
does not have factors $t_{j^+}^{+1}$. By Lemma \ref{lem1} (2) and Lemma \ref{jpm}, the exponent of $t_{j^+}$ is $0$ in $M$.
Since $M\in V(DG_l)\subset \mathcal{M}$, there is
a corresponding $\textbf{i}$-trail 
$\pi=(\gamma_0,\gamma_1,\cdots,\gamma_N)$. We can take integers $c_m$ ($m=1,2,\cdots,N$)
as
\[
\gamma_{m-1}-\gamma_m=c_m\alpha_{i_m}.
\]
By Proposition \ref{term-prop}, it holds
\[
c_j=0,\ c_{j^+}=1,\ \text{for }s\in[j+1,j^+-1]\ \text{such that }a_{i_s,i_j}<0,\ \text{one get }c_s=0.
\]
Let $\pi'=(\gamma_0',\gamma_1',\cdots,\gamma_N')$ be
the pre-\textbf{i}-trail and $c_m'$ ($m=1,2,\cdots,N$) be integers satisfying
\[
\gamma_{m-1}'-\gamma_m'=c_m'\alpha_{i_m}
\]
and 
\[
c_m'=c_m \quad \text{for}\ m\in[1,N]\setminus \{j,j^+\},
\]
\[
c_j'=c_j+1 \quad c_{j^+}'=c_{j^+}-1.
\]
Since $e_te_{t'}=e_{t'}e_t$ if $a_{t,t'}=0$, we can verify that $\pi'$ is an \textbf{i}-trail.
By Proposition \ref{trail-prop}, the monomial corresponding to $\pi'$ is
\[
\mathcal{M}\ni t_1^{d_1(\pi')}\cdots t_N^{d_N(\pi')}=t_1^{d_1(\pi)}\cdots t_N^{d_N(\pi)}A^{-1}_{j}=M\cdot A_j^{-1}.
\]
By the above argument, it holds $V(DG_{l+1})\subset \mathcal{M}$ so that $V(DG)=V(DG_r)\subset \mathcal{M}$. Note
that in step (2) of Theorem \ref{thm1}, if a monomial $M$ other than $t_{J} t_{J+1}^{a_{i_{J+1},i}}\cdots t_{N}^{a_{i_N,i}}$ in step (3)
is a sink of $DG_l$ then at least one arrow $M\rightarrow M'$ is added by Proposition \ref{term-prop2} and \ref{term-prop2-a}.

Next, we prove $V(DG)\supset \mathcal{M}$. For any $M=t_1^{d_1}\cdots t_N^{d_N}\in\mathcal{M}$, let us prove 
$M\in V(DG)$. 
Our claim is evident when $M=t_k$, so we assume that $M\neq t_k$.
Using Lemma \ref{lem3}, one can take a monomial $M_1\in \mathcal{M}$
and $m_1\in [1,N]$ such that $m_1^+\leq N$, $M=M_1\cdot A_{m_1}^{-1}$, $L(M_1)>L(M)$
and the exponents of $t_{m_1}$ and $t_{m_1^+}$ in $M_1$ are $1$ and $0$, respectively,
where $L$ is the notation in (\ref{L-def}).
Using Lemma \ref{lem3} repeatedly, since the set $\mathcal{M}$ is finite,
we can also take monomials
$M_2$, $M_3,\cdots,M_s=t_k\in \mathcal{M}$ and $m_2,\cdots,m_s$ such that
\[
M_{l-1}=M_{l}\cdot A_{m_l}^{-1},
\]
$m_l^+\leq N$, $L(M_{l})>L(M_{l-1})$ and the exponents of $t_{m_l}$ and $t_{m_l^+}$ in $M_l$ are $1$ and $0$, respectively
for $l=2,3,\cdots,s$. Since the exponent of $t_{m_l}$ is positive and $m_l^+ \leq N$,
we see that $M_l\neq t_{J} t_{J+1}^{a_{i_{J+1},i}}\cdots t_{N}^{a_{i_N,i}}$ for $l=1,2,\cdots,s$.
Thus, $M=t_k\cdot A_{m_s}^{-1}A_{m_{s-1}}^{-1}\cdots A_{m_2}^{-1}A_{m_1}^{-1} \in V(DG)$. \qed

\subsection{Properties of the decoration graph $DG$}

In this subsection, we prove that the decoration graph
$DG$ in Theorem \ref{thm1} possesses several similar properties to crystal graphs $B(\Lambda_i)$ for minuscule weights $\Lambda_i$.
First, the following proposition is an easy consequence of Theorem \ref{thm1}.

\begin{prop}
The directed graph $DG$ has the unique source $t_k$ and unique sink $t_{J} t_{J+1}^{a_{i_{J+1},i}}\cdots t_{N}^{a_{i_N,i}}$.
\end{prop}

Next, we prove the following:

\begin{prop}
In the decoration graph $DG$ of Theorem \ref{thm1}, we suppose that
there are arrows
$M \rightarrow M\cdot A_j^{-1}$ and $M \rightarrow M\cdot A_m^{-1}$ with $j,m\in[1,N]$ such that $j\neq m$ and $j^+\leq N$, $m^+\leq N$:
\[
\begin{xy}
(0,0) *{M}="1",
(10,-10)*{M\cdot A_j^{-1}}="2",
(-10,-10)*{M\cdot A_m^{-1}}="3",
\ar@{->} "1";"2"
\ar@{->} "1";"3"
\end{xy}
\]
Then there are two arrows $M\cdot A_m^{-1} \rightarrow M\cdot A_m^{-1}A_j^{-1}$ and $M\cdot A_j^{-1} \rightarrow M\cdot A_j^{-1}A_m^{-1}$
in $DG$:
\[
\begin{xy}
(0,0) *{M}="1",
(10,-10)*{M\cdot A_j^{-1}}="2",
(-10,-10)*{M\cdot A_m^{-1}}="3",
(0,-20)*{M\cdot A_j^{-1}A_m^{-1}}="4",
\ar@{->} "1";"2"
\ar@{->} "1";"3"
\ar@{->} "2";"4"
\ar@{->} "3";"4"
\end{xy}
\]
\end{prop}

\nd
{\it Proof.} We may assume $m<j$. Let $M=\prod^N_{l=1} t_l^{d_l}$ and $\pi=(\gamma_0,\cdots,\gamma_N)$ be the corresponding $\textbf{i}$-trail with nonnegative integers $\{c_l\}_{l\in[1,N]}$ satisfying $\gamma_{l-1}-\gamma_l=c_l\alpha_{i_l}$. By the assumption, Lemma \ref{lem1}, Lemma \ref{jpm} and (2) of Theorem \ref{thm1}, we see that
$d_j=d_m=1$, $d_{j^+}=d_{m^+}=0$. If $m^+<j$ then the exponents of $t_j$ and $t_{j^+}$ in $M\cdot A_{m}^{-1}$ are $1$ and $0$, respectively.
Thus, there is an arrow $M\cdot A_{m}^{-1}\rightarrow M\cdot A_{m}^{-1}A_j^{-1}$ by (2) of Theorem \ref{thm1}. Similarly,
there is an arrow $M\cdot A_{j}^{-1}\rightarrow M\cdot A_{j}^{-1}A_m^{-1}$. Hence, we assume $j<m^+$ so that $i_j\neq i_m$.
It follows from Proposition \ref{term-prop} that
\begin{equation}\label{DG-pr-1}
c_j=0,\ c_{j^+}=1,\ \text{for }s\in[j+1,j^+-1]\ \text{such that }a_{i_s,i_j}<0,\ \text{it holds }c_s=0,
\end{equation}
\begin{equation}\label{DG-pr-2}
c_m=0,\ c_{m^+}=1,\ \text{for }s\in[m+1,m^+-1]\ \text{such that }a_{i_s,i_m}<0,\ \text{it holds }c_s=0.
\end{equation}
In particular, it holds $c_{j^+}=c_{m^+}=1$.
We see that $a_{i_j,i_m}=0$ because if not, $m^+<j^+$ contradicts (\ref{DG-pr-1}) and $j^+<m^+$ contradicts (\ref{DG-pr-2}).
Therefore, the exponents of $t_j$, $t_{j^+}$ in $M$ are same as in $M\cdot A_m^{-1}$, which yields
there is an arrow $M\cdot A_{m}^{-1}\rightarrow M\cdot A_{m}^{-1}A_j^{-1}$ by (2) of Theorem \ref{thm1}.
Similarly,
there is an arrow $M\cdot A_{j}^{-1}\rightarrow M\cdot A_{j}^{-1}A_m^{-1}$. \qed

\vspace{2mm}

The above proposition is an analog of the following property of crystal graphs $B(\Lambda_i)$ for minuscule weights $\Lambda_i$:
If $b\in B(\Lambda_i)$ satisfies $\tilde{f}_lb\neq0$ and $\tilde{f}_rb\neq0$ with some $l,r\in I$ ($l\neq r$):
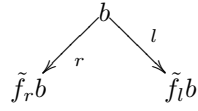
\[
\begin{xy}
(0,0) *{b}="1",
(10,-10)*{\tilde{f}_lb}="2",
(-10,-10)*{\tilde{f}_rb}="3",
\ar@{->} "1";"2"^l
\ar@{->} "1";"3"^r
\end{xy}
\]
then it holds 
$\tilde{f}_l\tilde{f}_rb=\tilde{f}_r\tilde{f}_lb\in B(\Lambda_i)$:
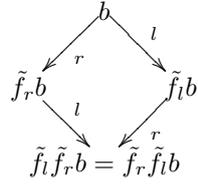
\[
\begin{xy}
(0,0) *{b}="1",
(10,-10)*{\tilde{f}_lb}="2",
(-10,-10)*{\tilde{f}_rb}="3",
(0,-20)*{\tilde{f}_l\tilde{f}_rb=\tilde{f}_r\tilde{f}_lb}="4",
\ar@{->} "1";"2"^l
\ar@{->} "1";"3"^r
\ar@{->} "2";"4"^r
\ar@{->} "3";"4"^l
\end{xy}
\]
Note that arrows in $DG$ are naturally colored by $[1,N]$ and arrows in
crystal graphs are colored by $I$.
As seen in the graph of $i=1$ of Example \ref{D4ex}, in $DG$
it may happen there are arrows $M \rightarrow M\cdot A_j^{-1}$ and $M \rightarrow M\cdot A_m^{-1}$
with $j,m\in[1,N]$ such that $i_j=i_m \in I$. In fact, there are two arrows
colored by $2$ and $9$ ($i_2=i_9=2$) starting from the monomial
$M=\frac{t_2t_9}{t_{10}t_{11}}$. 

There are examples the graph $DG$ coincides
with a subgraph of a genuine crystal graph of the monomial realization (see \ref{cry-mono}) for a minuscule representation
if we replace each color $l$ of arrow in $DG$ with $i_l$.
For instance, let $\mathfrak{g}$ be of type ${\rm A}_4$ and $\textbf{i}=(2,1,3,2,4,3,4,1,2,1)$. Then
the graph $DG$ for $\Delta_{w_0\Lm_3,s_3\Lm_3}\circ \theta^-_{\textbf{i}}(t_1,\cdots,t_{10})$ is as follows:
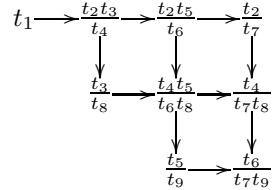
\[
\begin{xy}
(0,0) *{t_1}="1",
(10,0)*{\frac{t_2t_3}{t_4}}="2",
(20,0)*{\frac{t_2t_5}{t_6}}="3",
(10,-10)*{\frac{t_3}{t_8}}="4",
(20,-10)*{\frac{t_4t_5}{t_6t_8}}="5",
(30,0)*{\frac{t_2}{t_7}}="6",
(30,-10)*{\frac{t_4}{t_7t_8}}="7",
(20,-20)*{\frac{t_5}{t_9}}="8",
(30,-20)*{\frac{t_6}{t_7t_9}}="9",
\ar@{->} "1";"2"
\ar@{->} "2";"3"
\ar@{->} "2";"4"
\ar@{->} "3";"5"
\ar@{->} "4";"5"
\ar@{->} "3";"6"
\ar@{->} "6";"7"
\ar@{->} "5";"7"
\ar@{->} "5";"8"
\ar@{->} "8";"9"
\ar@{->} "7";"9"
\end{xy}
\]
Replacing variables as
$(t_1,t_2,\cdots,t_{10})\leftrightarrow (Y_{1,2},Y_{1,1},Y_{1,3},Y_{2,2},Y_{1,4},Y_{2,3},Y_{2,4},Y_{2,1},Y_{3,2},Y_{3,1})$
and coloring each arrow as $M\overset{i_j}{\rightarrow} M\cdot A_{j}^{-1}$,
the graph $DG$ coincides with the subgraph of monomial realization of $B(\Lambda_{2})$
with $p_{2,1}=p_{2,3}=p_{3,4}=1$ and the highest weight vector $Y_{1,2}$,
which is obtained by removing the lowest weight vector
$\frac{1}{Y_{3,3}}$ from the monomial realization.

\section{Examples of non-minuscule cases}\label{ex6-1}

\subsection{$\Delta_{w_0\Lm_i,s_i\Lm_i}\circ \theta^-_{\textbf{i}}(t_1,\cdots,t_N)$ of type ${\rm G}_2$}

In this section, 
we assume $\mathfrak{g}$ is of type ${\rm G}_2$ with Cartan matrix $(a_{i,j})_{i,j\in I}$ such that $a_{2,1}=-3$
and $a_{1,2}=-1$. 
Let us prove that the algorithm in Theorem \ref{thm1} works for this case.

\vspace{3mm}

\nd
\underline{In the case $\textbf{i}=(1,2,1,2,1,2)$}

\vspace{2mm}

By the results in \cite{KN},
\begin{eqnarray*}
& &\Delta_{w_0\Lm_1,s_1\Lm_1}\circ \theta^-_{\textbf{i}}(t_1,\cdots,t_6)\\
&=&t_1+\frac{t_2^3}{t_3}+3\frac{t_2^2}{t_4}+3\frac{t_2t_3}{t_4^2}+\frac{t_3^2}{t_4^3}
+2\frac{t_3}{t_5}+\frac{t_4^3}{t_5^2}+3\frac{t_2t_4}{t_5}+3\frac{t_2}{t_6}
+3\frac{t_3}{t_4t_6}+3\frac{t_4^2}{t_5t_6}+3\frac{t_4}{t_6^2}+\frac{t_5}{t_6^3}.
\end{eqnarray*}
\[
\Delta_{w_0\Lm_2,s_2\Lm_2}\circ \theta^-_{\textbf{i}}(t_1,\cdots,t_6)=t_6.
\]
Applying our algorithm in Theorem \ref{thm1} to $\Delta_{w_0\Lm_1,s_1\Lm_1}$, we get the following graph $DG$:
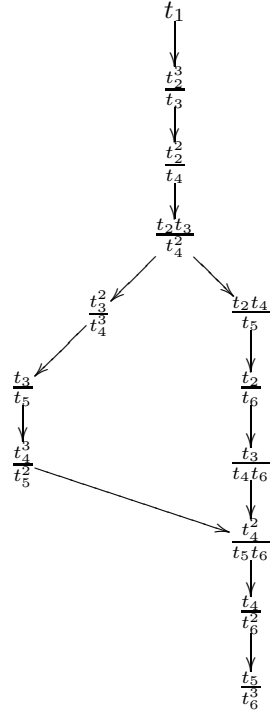
\[
\begin{xy}
(100,100) *{t_1}="1",
(100,90)*{\frac{t_2^3}{t_3}}="2",
(100,80)*{\frac{t_2^2}{t_4}}="3",
(100,70)*{\frac{t_2t_3}{t_4^2}}="4",
(90,60)*{\frac{t_3^2}{t_4^3}}="4-1",
(110,60)*{\frac{t_2t_4}{t_5}}="4-2",
(80,50)*{\frac{t_3}{t_5}}="5-1",
(110,50)*{\frac{t_2}{t_6}}="5-2",
(80,40)*{\frac{t_4^3}{t_5^2}}="6-1",
(110,40)*{\frac{t_3}{t_4t_6}}="6-2",
(110,30)*{\frac{t_4^2}{t_5t_6}}="7",
(110,20)*{\frac{t_4}{t_6^2}}="8",
(110,10)*{\frac{t_5}{t_6^3}}="9",
\ar@{->} "1";"2"^{}
\ar@{->} "2";"3"^{}
\ar@{->} "3";"4"^{}
\ar@{->} "4";"4-1"^{}
\ar@{->} "4";"4-2"^{}
\ar@{->} "4-1";"5-1"^{}
\ar@{->} "4-2";"5-2"^{}
\ar@{->} "5-1";"6-1"^{}
\ar@{->} "5-2";"6-2"^{}
\ar@{->} "6-1";"7"^{}
\ar@{->} "6-2";"7"^{}
\ar@{->} "7";"8"^{}
\ar@{->} "8";"9"^{}
\end{xy}
\]
Thus, our algorithm works in this case.
Note that the above graph coincides with the subgraph of crystal graph of a monomial realization of
$B(\Lambda_1)$ with $p_{1,2}=1$ by replacing variables $(t_1,\cdots,t_6)$
as $(Y_{1,1},Y_{1,2},Y_{2,1},Y_{2,2},Y_{3,1},Y_{3,2})$ and coloring arrows properly.
The algorithm for $\Delta_{w_0\Lm_2,s_2\Lm_2}\circ \theta^-_{\textbf{i}}(t_1,\cdots,t_6)$ clearly works.

\vspace{3mm}

\nd
\underline{In the case $\textbf{i}=(2,1,2,1,2,1)$}

\vspace{2mm}

Rewriting $\theta^-_{2,1,2,1,2,1}(t_1,\cdots,t_6)$ as
\[
\theta^-_{\textbf{i}}(t_1,\cdots,t_6)=
\theta^-_{2,1,2,1,2,1}(t_1,\cdots,t_6)
=\theta^-_{1,2,1,2,1,2}(t_1',\cdots,t_6'),
\]
it holds
\[
\Delta_{w_0\Lm_2,s_2\Lm_2}\circ \theta^-_{\textbf{i}}(t_1,\cdots,t_6)
=\Delta_{w_0\Lm_2,s_2\Lm_2}\circ \theta^-_{1,2,1,2,1,2}(t_1',\cdots,t_6')
=t_6'.
\]
It follows by Proposition 6.2 in \cite{N0} that
\[
\Delta_{w_0\Lm_2,s_2\Lm_2}\circ \theta^-_{\textbf{i}}(t_1,\cdots,t_6)
=t_6'=
t_1+\frac{t_2}{t_3}+\frac{t_3^2}{t_4}+2\frac{t_3}{t_5}+\frac{t_4}{t_5^2}+\frac{t_5}{t_6}.
\] 
Our algorithm to $\Delta_{w_0\Lm_2,s_2\Lm_2}\circ \theta^-_{\textbf{i}}(t_1,\cdots,t_6)$ computes the following graph $DG$:
\[
t_1\rightarrow
\frac{t_2}{t_3}\rightarrow
\frac{t_3^2}{t_4}\rightarrow
\frac{t_3}{t_5}\rightarrow
\frac{t_4}{t_5^2}\rightarrow
\frac{t_5}{t_6}.
\] 
Therefore, our algorithm works. We can check that above graph coincides with a subgraph of the crystal graph of
a monomial realization of $B(\Lambda_2)$ with $p_{2,1}=1$ by replacing variables $(t_1,\cdots,t_6)$
as $(Y_{1,2},Y_{1,1},Y_{2,2},Y_{2,1},Y_{3,2},Y_{3,1})$ and coloring arrows properly. Clearly, it holds
\[
\Delta_{w_0\Lm_1,s_1\Lm_1}\circ \theta^-_{\textbf{i}}(t_1,\cdots,t_6)=t_6
\]
and our algorithm works.

Note that coefficients greater than $1$ are appearing in contrary to the case of previous section.

\subsection{Specific reduced words}\label{ex6-2}

We found examples our algorithm works for a specific reduced word even if $V(-w_0\Lambda_i)$ is not minuscule.

\begin{ex}\label{non-mini-ex1}

Let $\mathfrak{g}$ be of type ${\rm B}_3$ and $\textbf{i}=(1,2,3,1,2,3,1,2,3)$. Note that $V(-w_0\Lambda_1)$ and $V(-w_0\Lambda_2)$ are
not minuscule, however,
applying the algorithm in Theorem \ref{thm1} to $i=1,2$, one obtain the following:
\[
t_1\rightarrow \frac{t_2}{t_4}\rightarrow \frac{t_3^2}{t_5}\rightarrow \frac{t_3}{t_6}
\rightarrow \frac{t_5}{t_6^2}\rightarrow \frac{t_7}{t_8}\qquad (i=1),
\]
\[
t_4\rightarrow \frac{t_5}{t_7}\rightarrow \frac{t_6^2}{t_8}\rightarrow \frac{t_6}{t_9}
\rightarrow \frac{t_8}{t_9^2}\qquad (i=2).
\]
By the results in [Theorem 8.1, \cite{N1}], we obtain 
\[
\Delta_{w_0\Lambda_1,s_1\Lambda_1}\circ \theta^-_{\textbf{i}}(t_1,t_2,\cdots,t_9)=
t_1+\frac{t_2}{t_4}+\frac{t_3^2}{t_5}+2\frac{t_3}{t_6}+\frac{t_5}{t_6^2}+\frac{t_7}{t_8},
\]
\[
\Delta_{w_0\Lambda_2,s_2\Lambda_2}\circ \theta^-_{\textbf{i}}(t_1,t_2,\cdots,t_9)=
t_4+\frac{t_5}{t_7}+\frac{t_6^2}{t_8}+2\frac{t_6}{t_9}+\frac{t_8}{t_9^2}.
\]
Hence, for $i=1,2$,
the set of monomials appearing in
$\Delta_{w_0\Lm_i,s_i\Lm_i}\circ \theta^-_{\textbf{i}}(t_1,\cdots,t_N)$
coincides with the set of Laurent monomials in the directed graph $DG$,
which implies the algorithm in Theorem \ref{thm1} works in this setting.

\end{ex}

\begin{ex}\label{non-mini-ex2}

Let $\mathfrak{g}$ be of type ${\rm C}_3$ and $\textbf{i}=(2,1,3,2,1,3,2,1,3)$. Note that $V(-w_0\Lambda_2)$ is
not minuscule, however,
applying the algorithm in Theorem \ref{thm1} to $i=2$, one obtain the following:
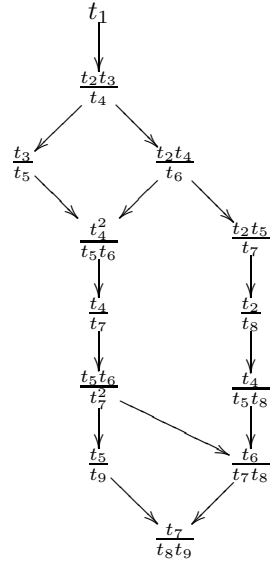
\[
\begin{xy}
(0,100) *{t_1}="1",
(0,90)*{\frac{t_2t_3}{t_4}}="2",
(-10,80)*{\frac{t_3}{t_5}}="2-1",
(10,80)*{\frac{t_2t_4}{t_6}}="2-2",
(0,70)*{\frac{t_4^2}{t_5t_6}}="3-1",
(20,70)*{\frac{t_2t_5}{t_7}}="3-2",
(0,60)*{\frac{t_4}{t_7}}="4-1",
(20,60)*{\frac{t_2}{t_8}}="4-2",
(0,50)*{\frac{t_5t_6}{t_7^2}}="5-1",
(20,50)*{\frac{t_4}{t_5t_8}}="5-2",
(0,40)*{\frac{t_5}{t_9}}="6-1",
(20,40)*{\frac{t_6}{t_7t_8}}="6-2",
(10,30)*{\frac{t_7}{t_8t_9}}="7",
\ar@{->} "1";"2"^{}
\ar@{->} "2";"2-1"^{}
\ar@{->} "2";"2-2"^{}
\ar@{->} "2-1";"3-1"^{}
\ar@{->} "2-2";"3-1"^{}
\ar@{->} "2-2";"3-2"^{}
\ar@{->} "3-1";"4-1"^{}
\ar@{->} "3-2";"4-2"^{}
\ar@{->} "4-1";"5-1"^{}
\ar@{->} "4-2";"5-2"^{}
\ar@{->} "5-1";"6-1"^{}
\ar@{->} "5-2";"6-2"^{}
\ar@{->} "5-1";"6-2"^{}
\ar@{->} "6-1";"7"^{}
\ar@{->} "6-2";"7"^{}
\end{xy}
\]
On the other hands, by using the bilinear form (\ref{minor-bilin}), one can compute
\begin{eqnarray*}
\Delta_{w_0\Lambda_2,s_2\Lambda_2}\circ \theta^-_{\textbf{i}}(t_1,t_2,\cdots,t_9)
&=&t_1+\frac{t_2t_3}{t_4}+\frac{t_3}{t_5}+\frac{t_2t_4}{t_6}+\frac{t_4^2}{t_5t_6}+\frac{t_2t_5}{t_7}\\
& &+2\frac{t_4}{t_7}+\frac{t_2}{t_8}+\frac{t_5t_6}{t_7^2}+\frac{t_4}{t_5t_8}+\frac{t_5}{t_9}
+\frac{t_6}{t_7t_8}+\frac{t_7}{t_8t_9}.
\end{eqnarray*}
Hence, our algorithm works.
\end{ex}

\begin{rem}\label{last-rem}

In the case $V(\Lambda_i)$ is not minuscule, the algorithm in Theorem \ref{thm1} does not work in general, but, 
there are several examples it works such as type ${\rm G_2}$ case or Example \ref{non-mini-ex1}, \ref{non-mini-ex2}.
For each $j,l\in I$ such that
$a_{j,l}<0$, if the subword of $\textbf{i}=(i_1,\cdots,i_N)$ consisting of all $j$, $l$ is $(\cdots,j,l,j,l,j,l)$ or $(\cdots,l,j,l,j,l,j)$
then we expect the algorithm
computes all monomials in $\Delta_{w_0\Lambda_i,s_i\Lambda_i}\circ \theta^-_{\textbf{i}}(t_1,t_2,\cdots,t_N)$ correctly.
Furthermore, we also expect the all monomials in $DG$ appear in the
crystal graph of a monomial realization for $B(\Lambda_{i_k})$
with
\[
p_{j,l}=
\begin{cases}
1 & \text{if the subword is } (j,l,j,l,j,l,\cdots), \\
0  & \text{if the subword is } (l,j,l,j,l,j,\cdots)
\end{cases}
\]
by replacing variables $t_m\leftrightarrow Y_{s,r}$ if $i_m=r$ and $s=\sharp\{a\in[1,m] | i_a = r\}$.
In fact,
the reduced words in type ${\rm G}_2$ case, Example \ref{non-mini-ex1}, \ref{non-mini-ex2} and
the example in the end of Sect. \ref{minus} satisfy
the above assumption and the algorithm works and monomials in $DG$ coincide with
monomials in a subgraph of crystal graph of a monomial realization.
\end{rem}

\end{document}